\documentclass{amsart}
\usepackage{amsmath,amssymb}
\usepackage[all]{xy}

\theoremstyle{plain}
\newtheorem{theorem}{Theorem}[section]
\newtheorem{proposition}[theorem]{Proposition}
\newtheorem{lemma}[theorem]{Lemma}
\newtheorem{corollary}[theorem]{Corollary}

\newtheorem{problem}[theorem]{Problem}

\newtheorem{theorem*}{Theorem}[]
\theoremstyle{definition}
\newtheorem{definition}[theorem]{Definition}
\newtheorem{example}[theorem]{Example}
\newtheorem{examples}[theorem]{Examples}

\theoremstyle{remark}

\newcommand{\secref}[1]{Section~\ref{#1}}
\newcommand{\thmref}[1]{Theorem~\ref{#1}}
\newcommand{\propref}[1]{Proposition~\ref{#1}}
\newcommand{\lemref}[1]{Lemma~\ref{#1}}
\newcommand{\corref}[1]{Corollary~\ref{#1}}

\newcommand{\exref}[1]{Example~\ref{#1}}

\newcommand{\defref}[1]{Definition~\ref{#1}}

\def\R{{\mathbb \R}}

\def\Q{{\mathbb Q}}

\def\map{\mathrm{Map}}
\def\top{\mathrm{Top}}
\def\diff{\mathrm{Diff}}

\def\Hom{\mathrm{Hom}}

\def\cat0{\mathrm{cat}_0}

\def\ker{\mathrm{ker}\,}

\def\im{\mathrm{im}\,}
\def\EM{\text{Eilenberg-Mac$\:$Lane }}
\def\ev{\mathrm{w}}

\begin{document}

\title[Evaluation Maps in Rational Homotopy]{Evaluation Maps in Rational Homotopy}

\author{Yves F{\'e}lix}

\address{Institut Math{\'e}matique,
Universit{\'e} Catholique de Louvain, B-1348 Louvain-la-Neuve
Belgique}

\email{felix@math.ucl.ac.be}

\author{Gregory  Lupton}

\address{Department of Mathematics,
          Cleveland State University,
          Cleveland OH 44115 U.S.A.}

\email{G.Lupton@csuohio.edu}

\date{\today}



\begin{abstract}
Let $E$ be an $H$-space acting on a based space $X$. Then we refer
to $\ev \colon E \to X$, the map obtained by acting on the base
point of $X$, as a ``generalized evaluation map"  (see
\defref{def:evaluation map} for a precise definition).
We establish several fundamental results about the rational homotopy
behaviour of a generalized evaluation map, all of which apply to the
usual evaluation map $\map(X, X;1) \to X$. With mild hypotheses on
$X$, we show that a generalized evaluation map $\ev$ factors, up to
rational homotopy, through a map $\Gamma_\ev \colon S_\ev \to X$
where $S_\ev$ is a (relatively small) finite product of
odd-dimensional spheres and $\pi_\#(\Gamma_\ev)\otimes \mathbb Q$ is
injective. This result has strong consequences: if the image in
rational homotopy groups of $\ev$ is trivial, then \emph{the
generalized evaluation map is null-homotopic} after rationalization;
unless $X$ satisfies a very strong splitting condition, any
generalized evaluation map \emph{induces the trivial homomorphism in
rational cohomology}; the map $\Gamma_\ev$ is \emph{rationally a
homotopy monomorphism} and a generalized evaluation map may be
written as a composition of a homotopy epimorphism and this homotopy
monomorphism. We include illustrative examples and prove numerous
subsidiary results of interest.
\end{abstract}

\maketitle

\section{Introduction}\label{sec:intro}

Let $X$ be a based space and let $\map(X, X)$ be the space of
unbased, or free, maps from $X$ to itself.  In general $\map(X,
X)$ is disconnected; we denote by $\map(X, X;1)$ its identity
component, that is, the path component that consists of self maps
that are (freely) homotopic to the identity.  Then we have
\emph{the evaluation map} $\omega\colon \map(X, X;1) \to X$
defined by evaluation at the basepoint of $X$. This map occupies a
central place in the homotopy theory of fibrations
(cf.~\cite{Got1,Got2,Go1,Go2}).

The evaluation map $\omega$ and its rationalization will play a
distinguished role in this paper.  However, we find that our
methods and results apply equally well to other contexts in which
one has an ``evaluation map."  For example, it is often of
interest to consider the space $\top(X,X)$ of self-homeomorphisms
of $X$ and the corresponding evaluation map $\ev\colon \top(X,X;1)
\to X$. Here, $\top(X,X;1)$ denotes the component of $\top(X,X)$
that consists of self-homeomorphisms homotopic (\emph{via}
self-homeomorphisms) to the identity. Likewise, if $X$ is a smooth
manifold, then one may replace $\top(X,X)$ with $\diff(X,X)$, and
so-forth. A further example of an ``evaluation map" to which our
methods apply concerns \emph{configuration spaces}. Let $F(X,k)$
denote the configuration space that consists of ordered $k$-tuples
of distinct points in a space $X$, and let $(p_1, \dots, p_k)$ be
a choice of basepoint in $F(X,k)$. Then we have a map
$$\theta \colon \top(X,X;1) \to F(X,k)$$
given by $\theta(\alpha) = (\alpha(p_1), \dots, \alpha(p_k))$.
Actually, here we have $\theta = \ev\circ\Theta$, where
$$\Theta \colon \top(X,X;1) \to \top\big(F(X,k),F(X,k); 1\big)$$
is the natural injection defined by $\Theta(\alpha)(q_1, \dots,
q_k) = (\alpha(q_1), \dots, \alpha(q_k))$, and
$$\ev \colon \top\big(F(X,k), F(X,k); 1\big) \to F(X,k)$$
is an evaluation map for $F(X,k)$ in the preceding sense, with
``$\top$" replacing ``$\map$."

Motivated by the preceding examples, we now make a formal
definition of the evaluation maps that we consider.  Recall that
an $H$-space is a pair $(E, \mu)$ with $E$ a based space and
multiplication $\mu \colon E \times E \to E$ a based map that
satisfies $\mu\circ J \sim \nabla \colon E \vee E \to E$. Here,
$\nabla \colon E \vee E \to E$ denotes the folding map and $J
\colon E \vee E \to E \times E$ the obvious inclusion. We say that
the multiplication has \emph{strict identity} if $\mu\circ J =
\nabla$ (equals, not just homotopic).  Note that $\map(X,X;1)$ is
an $H$-space with strict identity.   Now let $i_1\colon E \to
E\times X$ and $i_2\colon X \to E\times X$ denote the inclusions.
By an \emph{action of $E$ on $X$} we mean a map $A\colon E\times X
\to X$ that satisfies $A\circ i_2 = 1\colon X \to X$.  We say that
\emph{the action is associative} if in addition we have
$A\circ(\mu\times1) = A\circ(1\times A)$.

\begin{definition}\label{def:evaluation map}
A \emph{generalized evaluation map} is any (based) map $\ev \colon
E \to X$, from a connected $H$-space with strict identity $E$ to a
space $X$, for which there exists an associative action $A\colon
E\times X \to X$ that restricts to $\ev$, that is, that satisfies
$A\circ i_1 = \ev \colon E \to X$.
\end{definition}

\begin{examples}
(1) The action $A \colon \map(X,X;1) \times X \to X$ given by
$A(f, x) = f(x)$ makes $\omega\colon \map(X,X;1) \to X$ a
generalized evaluation map according to \defref{def:evaluation
map}. Similarly for all the other examples mentioned above.

(2) Suppose $G$ is a connected topological group and $A\colon G
\times X \to X$ is a group action in the usual sense.  Then the
\emph{orbit map} of the action is a generalized evaluation map $G
\to X$.

(3) More generally, suppose given a fibration $X \to Y \to B$.
Then the connecting map $\partial\colon \Omega B \to X$ is a
generalized evaluation map.  This follows from the usual action of
$\Omega B$ on the fibre $X$.  Note, however, that we must take
Moore loops in $\Omega X$ to obtain an $H$-space with strict
identity.
\end{examples}

Revert now to the ordinary evaluation map $\omega \colon \map(X,X;1)
\to X$.  \emph{For the remainder of the paper, we assume that $X$ is a
nilpotent, finite complex}. Since $X$ is finite, a result of Milnor
\cite{Mil} implies that $\map(X, X;1)$ is a CW complex. Since $X$ is
nilpotent, we may choose and fix a rationalization $e\colon X \to
X_\Q$.  Now results of \cite{H-M-R} imply that $e_*\colon \map(X,X;1)
\to \map(X,X_\Q;e)$ is a rationalization. Thus, the map $\omega_\Q
\colon \map(X,X_\Q;e) \to X_\Q$, also defined by evaluation at the
basepoint of $X$, may be taken to be the rationalization of $\omega$.
We refer to $\omega_\Q$ as \emph{the rationalized evaluation map}.
Recall that the \emph{$n$th Gottlieb group of $X$}, denoted $G_n(X)$,
is the subgroup of $\pi_n(X)$ defined as the image of $\omega_\# \colon
\pi_n\big(\map(X, X; 1)\big) \to \pi_n(X)$ \cite{Go1}. The subgroup of
$\pi_n(X_\Q)$ defined as the image of $(\omega_\Q)_\# \colon
\pi_n\big(\map(X, X_\Q; e)\big) \to \pi_n(X_\Q)$ is called the
\emph{$n$th rationalized Gottlieb group of $X$} and denoted by
$G_n(X_\Q)$. By a theorem of Lang \cite{Lan}, we have $G_n(X_\Q) \cong
G_n(X)\otimes\Q$ under our assumption that $X$ is finite. The
rationalized Gottlieb groups have played an important role in some of
the major developments of rational homotopy theory (cf.~\cite{F-H,
Hal2}). Our results in this paper show that the rationalized Gottlieb
groups exercise a very strong determining effect on the rationalized
evaluation map.

A result of F{\'e}lix-Halperin (\cite[Th.III]{F-H}) implies that
$G_{2i}(X_\Q) = 0$ for all $i$ and $G_{2i+1}(X_\Q)$ is non-zero for
only finitely many $i$. Suppose $\{\alpha_1, \alpha_2, \ldots,
\alpha_r\}$ is a basis of $G_{*}(X_\Q) = G_{\text{odd}}(X_\Q)$ with
$\alpha_i \in G_{n_i}(X_\Q)$ (here we regard an element of
$\pi_n(X_\Q)$ as represented by a map $\alpha \colon S^n_\Q \to X_\Q$).
For each $\alpha_i$, we may choose a $\beta_i \in
\pi_{n_i}\big(\map(X,X_\Q;e)\big)$ such that $\omega_\Q\circ\beta_i =
\alpha_i$.  The adjoint of $\beta_i$ gives a map $F_i \colon S_\Q^{n_i}
\times X \to X_\Q$ that extends the map $(\alpha_i \mid e) \colon
S_\Q^{n_i} \vee X \to X_\Q$. Denote by $S_X$ the product of
odd-dimensional rational spheres $S_\Q^{n_1} \times \cdots \times
S_\Q^{n_r}$ whose factors correspond to the domains of the basis of
$G_{*}(X_\Q)$. Then we form a map $F \colon S_X \times X \to X_\Q$ as
the composition
$$F = F_1\circ(1\times F_2)\circ \cdots \circ (1\times \cdots
\times1\times F_r).$$
Now set $\Gamma_X = F\circ i \colon S_X \to X_\Q$, where $i$ denotes
the inclusion of the product of spheres as the first $r$ factors.  We
refer to $\Gamma_X$ as a \emph{total Gottlieb element of $X_\Q$}. By
taking the adjoint of $F$, we obtain a lift $\widetilde{\Gamma}_X\colon
S_X \to \map(X,X_\Q;e)$ of $\Gamma_X$ through the rationalized
evaluation map $\omega_\Q$.

We prove the following result:

\begin{theorem}\label{thm:Gamma factor}
Let $X$ be any nilpotent, finite complex. The rationalized evaluation
map $\omega_{\Q}\colon \map(X,X_\Q;e) \to X_\Q$ factors up to homotopy
through the total Gottlieb element $\Gamma_X\colon S_X \to X_\Q$.  In
fact, there is a retraction $r$ of $\widetilde{\Gamma}_X$, that is, a
map $r \colon \map(X,X_\Q;e) \to S_X$ with $r\circ \widetilde{\Gamma}_X
= 1$, such that $\omega_{\mathbb Q}= \Gamma_X\circ r$.
\end{theorem}

This basic result has several strong consequences.  An immediate
one is the following striking illustration of the effect that the
homomorphism induced on rational homotopy groups has on the
rationalized evaluation map.

\begin{corollary}\label{cor:G trivial}
The evaluation map $\omega \colon \map(X,X;1) \to X$ is rationally
null-homotopic if and only if $G_*(X_\Q) = 0$.
\end{corollary}

Now the evaluation map $\omega$ may be viewed as a ``universal
connecting map" for fibrations with fibre $X$, in that any
connecting map $\Omega B \to X$ of a fibration $X \to E \to B$
factors through $\omega$ \cite{Got2}.  A further immediate
consequence of \thmref{thm:Gamma factor}, therefore, is the
following result.

\begin{corollary}\label{cor:partial trivial}
Let $X \to E \to B$ be any fibration with fibre $X$ a nilpotent,
finite space.  If $G_*(X_\Q) = 0$, then the connecting map
$\partial\colon \Omega B \to X$ is rationally null-homotopic.
\end{corollary}

There are many spaces to which these corollaries may be applied.
For instance, any suspension that is not rationally equivalent to
a sphere has trivial rationalized Gottlieb groups. Roughly
speaking, a typical wedge or connected sum of spaces has trivial
Gottlieb groups, as do many non-elliptic, coformal spaces. More
precisely, a space whose rational homotopy Lie algebra has trivial
centre has trivial rationalized Gottlieb group. Therefore, by
\corref{cor:G trivial}, the rationalized evaluation map is
null-homotopic in all such cases.

The preceding discussion of $\omega$ and the Gottlieb groups
extends naturally to generalized evaluation maps. Suppose given
$\ev \colon E \to X$ any generalized evaluation map.  In
\secref{sec:w-factor} we construct a map $\Gamma_\ev \colon S_\ev
\to X_\Q$ such that $\im (\Gamma_{\ev})_\#\otimes \mathbb Q = \im
(\ev)_\#\otimes \mathbb Q$. As with $S_X$ above, $S_\ev$ is a
product of a relatively small number of odd-dimensional rational
spheres. We refer to $\Gamma_\ev$ as a \emph{total Gottlieb
element of $X_\Q$ with respect to $\ev$}. Furthermore,
$\Gamma_\ev$ admits a lift through $\ev_\Q$, the rationalization
of $\ev$.  That is, there exists a map $\widetilde{\Gamma}_\ev
\colon S_\ev \to E_{\Q}$ that satisfies $\ev_\Q\circ
\widetilde{\Gamma}_\ev = \Gamma_\ev$. Then we have generalizations
of \thmref{thm:Gamma factor} and \corref{cor:G trivial} as
follows.

\begin{theorem}\label{thm:Gamma+ factor}
Let $\ev \colon E \to X$ be any generalized evaluation map with
$X$ a nilpotent, finite complex.  Suppose that $\Gamma_\ev \colon
S_\ev \to X_\Q$ is a total Gottlieb element of $X_\Q$ with respect
to $\ev$. Then $\ev_{\Q}$ factors up to homotopy through
$\Gamma_\ev$. More precisely, suppose that
$\widetilde{\Gamma}_\ev\colon S_\ev \to E_\Q$ is a lift of
$\Gamma_\ev$ through $\ev_\Q$.  Then there is a retraction $r
\colon E_\Q \to S_\ev$ of $\widetilde{\Gamma}_\ev$ such that
$\ev_{\mathbb Q} = \Gamma_\ev\circ r$.
\end{theorem}

\begin{corollary}\label{cor:w-sharp0}
Let $\ev \colon E \to X$ be any generalized evaluation map. Then
$\ev_\#\otimes\Q = 0 \colon \pi_*(E_\Q) \to \pi_*(X_\Q)$ if and only if
$\ev \colon E \to X$ is rationally null-homotopic.
\end{corollary}

We continue with a theorem related to the homotopy behaviour of
the maps $\Gamma_\ev \colon S_\ev \to X_\Q$. Recall that a map $f
\colon X \to Y$ is a \emph{homotopy monomorphism} if, for any $A$,
the induced map of homotopy sets $f_* \colon [A,X] \to [A,Y]$ is
injective \cite{Gan1}. In general it is a difficult problem to
identify when a map is a homotopy monomorphism.  We say that a map
of nilpotent spaces $f\colon X \to Y$ is a \emph{homotopy
monomorphism in the nilpotent category} if $f_* \colon [A,X] \to
[A,Y]$ is injective whenever $A$ is a nilpotent space.

\begin{theorem}\label{thm:mono} Let $X$ be a nilpotent,
finite complex and $\ev \colon E \to X$ be any generalized
evaluation map. Then $\Gamma_\ev \colon S_\ev \to X_\Q$ is a
homotopy monomorphism in the nilpotent category.
\end{theorem}

\thmref{thm:mono} is proved towards the end of
\secref{sec:technical}. As a consequence, together with
\thmref{thm:Gamma+ factor} we find that, after rationalization, a
generalized evaluation map may be written as a composition $\ev_\Q =
\Gamma_\ev \circ r$ of a homotopy epimorphism and a homotopy
monomorphism in the nilpotent category (\corref{cor:epi-mono}). We
also note the following immediate consequence of \thmref{thm:mono}.

\begin{corollary}
Let $X$ be a nilpotent, finite complex and let $\alpha\colon S^n
\to X_\Q$ be any rationalized Gottlieb element.  Then $\alpha$ is
a homotopy monomorphism in the nilpotent category.
\end{corollary}

In particular, this implies that the rationalized Hopf maps are
homotopy monomorphisms in the nilpotent category.  By contrast, the
Hopf map $\eta\colon S^7 \to S^4$ is not a homotopy monomorphism
\cite{Gan1}.

A further consequence of \thmref{thm:mono} is the classification up
to rational homotopy of cyclic maps. A map $f \colon A \to X$ is
called \emph{cyclic} if $(f\mid 1) \colon A \vee X \to X$ extends to
a map $A \times X\to X$ \cite{Var}.    Denote by $G(A,X)$ the set of
homotopy classes of cyclic maps from $A$ into $X$. This is a
generalization of the $n$th Gottlieb group of $X$, which we obtain
by taking $A = S^n$. Upon rationalizing a cyclic map, we obtain a
map $f_\Q \colon A \to X_\Q$ in $G(A,X_\Q)$.

\begin{theorem}\label{thm:cyclic} Let $X$ be a nilpotent, finite complex and let $A$ be any nilpotent
space.  Then there is a bijection of sets
$$G(A,X_\Q) \cong [A,
S_X] \cong \oplus_r \mbox{\rm Hom} (H_r(A;\mathbb Q),
G_r(X_{\mathbb Q}))\,.$$
\end{theorem}

\noindent{}This classification allows us, for instance, to easily
identify situations in which $G(A,X_\Q)$ is trivial and, hence,
$G(A,X)$ is finite. In \thmref{thm:cyclic general}, we extend this
result to apply to any generalized evaluation map.

Our last topic is the (co)homological behaviour of generalized
evaluation maps.  For the ordinary evaluation map $\omega\colon
\map(X,X;1) \to X$, this behaviour has been studied by Gottlieb
\cite{Gott72} and Oprea \cite{Op86, Op87}. From \cite{Op86} we have
the following result:

\begin{theorem}[Oprea]\label{thm:Oprea}
Let $F \to E\to B$ be a fibration with connecting map $\partial
\colon \Omega B \to F$. Suppose that $B$ is $1$-connected and that
$B$ and $F$ have finite type rational homology. Then there is a
splitting, up to rational homotopy, $F \simeq_{\Q} S \times Y$
with $S$ a product of \EM spaces and $$\mbox{\rm dim}\,
\pi_*(S)\otimes \Q = \mbox{\rm dim Image}\, (h_F\circ  \partial_\#
\colon \pi_*(\Omega B)\otimes \Q \to H_*(F;\Q))\,.$$ Here,
$h_F\colon \pi_*(F)\otimes\Q \to H_*(F;\Q)$ denotes the rational
Hurewicz homomorphism.
\end{theorem}

Oprea's result may be applied to the evaluation map $\omega$ by
considering it as the connecting map in the universal fibration
for fibrations with fibre $X$. Our main result about the
homological behaviour of a generalized evaluation map is the
following composite theorem, which gives a complete description
for rational coefficients.

\begin{theorem}\label{th:fact.homology} Let $\ev \colon E \to X$
be any generalized evaluation map with $X$ a nilpotent, finite
complex.  Then we have:
\begin{enumerate}
\item $\widetilde{H}_*(\ev;\Q) \not=0
\colon \widetilde{H}_*(E;\Q) \to \widetilde{H}_*(X;\Q)$ if and only
if $h_X \circ (\ev_\#\otimes\Q) \not=0 \colon \pi_*(E)\otimes\Q \to
\pi_*(X)\otimes\Q \to H_*(X;\Q)$;
\item if $h_X \circ (\ev_\#\otimes\Q)$ has image in $H_*(X;\Q)$ of
dimension $r > 0$, then $H_*(\ev;\Q)$ has image in $H_*(X;\Q)$ of
dimension $2^r$ and there is a rational homotopy equivalence $X
\simeq_{\Q} S \times Y$, with $S$ a product of odd-dimensional
spheres such that $H_*(S;\Q) \cong \mbox{\rm Image}\, H_*(\ev;\Q)$
and $\pi_*(S)\otimes\Q \cong \mbox{\rm Image} \,h_X \circ
(\ev_\#\otimes\Q)$;
\item if $X \simeq_\Q S^{2n+1} \times Y$, then
$\widetilde{H}_*(\omega;\Q) \not=0$, where $\omega\colon
\map(X,X;1) \to X$ is the ordinary evaluation map.
\end{enumerate}
\end{theorem}

Our treatment here extends Oprea's theorem to a generalized
evaluation map. \thmref{th:fact.homology} shows that, in most
cases, the rank of $H_*(\ev;\Q)$ is relatively small. We also
deduce that $H_*(\ev;\Q)$ is surjective only when $X$ is an
$H_0$-space. \thmref{th:fact.homology} has various interesting
corollaries, such as the following sharpening of a result of
Gottlieb \cite[Th.3]{Gott72} for rational coefficients.

\begin{corollary}\label{cor:Gottlieb}
Suppose that $\chi(X) \not = 0$. Then for every generalized
evaluation map $\ev \colon E \to X$, we have
$\widetilde{H}_*(\ev;\Q)) = 0 \colon \widetilde{H}_*(E;\Q) \to
\widetilde{H}_*(X;\Q)$.
\end{corollary}

A further consequence is the following result:

\begin{corollary}\label{cor:symplectic}
Let $M$ be a simply connected, symplectic manifold.  Then every
generalized evaluation map $\ev \colon E \to M$ is trivial on
rational homology, that is, $\widetilde{H}_*(\ev;\Q) = 0 \colon
\widetilde{H}_*(E;\Q) \to \widetilde{H}_*(M;\Q)$. Consequently, if
$G$ is a connected Lie group and $a\colon G \to M$ is the orbit map
of any $G$-action on $M$, we have $\widetilde{H}_*(a;\Q) = 0$.
\end{corollary}

These corollaries appear as \corref{cor:Gottlieb general} and
\corref{cor:symplectic general}, respectively.

The text is divided into five parts. In \secref{sec:w-factor} we
present the factorization results. \secref{sec:technical} contains
some technical lemmas on Gottlieb groups, and the monomorphism
theorem. The homological behaviour of generalized evaluation maps
is discussed in \secref{sec:cohomology}. \secref{sec:examples} is
a brief, concluding section in which we mention several problems
that suggest directions for future work.

We finish this introduction with some terminology and notation. We
work in the homotopy category, and so we often do not distinguish
between a map and the homotopy class it represents.  We use
$\simeq$ to denote that two spaces are homotopy equivalent, or
that a map is a homotopy equivalence. If $f \colon A \to B$ is a
map, then $f^*$ denotes pre-composition by $f$ and $f_*$ denotes
post-composition by $f$. We use $H_*(f)$ and $H^*(f)$ to denote
the map induced on homology, respectively cohomology, by the map
of spaces $f$, and $f_\#$ to denote the map induced on homotopy
groups.  Likewise, $\widetilde{H}_*(f)$ and $\widetilde{H}^*(f)$
denote reduced (co)homology.   We denote the rationalization of a
space $X$ by $X_\Q$ and of a map $f$ by $f_\Q$ (cf.~\cite{H-M-R}).
By an $H_0$-space, we mean a space whose rationalization is an
$H$-space. We say that maps $f, g \colon X \to Y$ are rationally
homotopic if their rationalizations are homotopic.  We denote this
relation either by $f =_\Q g\colon X \to Y$ or by $f_\Q =
g_\Q\colon X_\Q \to Y_\Q$. We reserve $\omega$ to denote the
evaluation map $\omega \colon \map(X,X;1) \to X$. Generalized
evaluation maps will be denoted with a generic $\ev$.  For the
remainder of the paper, we will usually drop the ``generalized"
and refer simply to \emph{an evaluation map}.

We assume familiarity with rational homotopy theory and use the
standard notation and terminology for minimal models as presented
in \cite{F-H-T}.  The basic facts that we use are as follows: Each
nilpotent space $X$ has a unique Sullivan minimal model
$(\mathcal{M}_{X}, d_X)$ in the category of commutative DG
(differential graded) algebras over $\Q$.  This DG algebra
$(\mathcal{M}_{X}, d_X)$ is of the form $\mathcal{M}_{X} = \land
V$, a free graded commutative algebra generated by a positively
graded vector space $V$ of finite type.  The differential $d_X$ is
decomposable, in that $d_X(V) \subseteq \land^{\geq 2} V$, and $V$
admits a basis $\{v_{\alpha}\}$ indexed by a well ordered set such
that $d_X(v_\alpha) \in \land (\{v_\beta\}_{\beta <\alpha})$.  A
fact that we use very frequently here is that an $H_0$-space has a
minimal model with zero differential. Each map $f\colon X \to Y$
also has a Sullivan minimal model which is a DG algebra map
$\mathcal{M}_f \colon \mathcal{M}_Y \to \mathcal{M}_X$. The
Sullivan minimal model is a complete rational homotopy invariant
for a space or a map.  If $f, g \colon X \to Y$ are maps of
\emph{rational spaces}, then $f$ and $g$ are homotopic if and only
if their Sullivan minimal models $\mathcal{M}_f$ and
$\mathcal{M}_g$ are homotopic in an algebraic sense. Rational
cohomology is readily retrieved from Sullivan minimal models:  We
have a natural isomorphism $H(\mathcal{M}_{X}, d_X) \cong
H^*(X;\Q)$ and this isomorphism identifies $H(\mathcal{M}_f)
\colon H(\mathcal{M}_Y) \to H(\mathcal{M}_X)$ with $H^*(f) \colon
H^*(Y; \Q) \to H^*(X; \Q)$. Rational homotopy groups are retrieved
as follows:  Let $Q(\mathcal{M}_X) \cong V$ be the (quotient)
module of indecomposables of $\mathcal{M}_X$. There is a natural
isomorphism $Q(\mathcal{M}_X) \cong \Hom(\pi_*(X),\Q)$, that
identifies $Q(\mathcal{M}_f) \colon Q(\mathcal{M}_Y) \to
Q(\mathcal{M}_X)$ with $(f_\#\otimes\Q)^* \colon \Hom(\pi_*(Y),\Q)
\to \Hom(\pi_*(X),\Q)$.

\bigskip

\noindent{\textsc{acknowledgements}.} It is a pleasure to thank
John Oprea for fruitful discussions on the general topics of this
paper. We also thank Sam Smith, whose work with the second-named
author in \cite{Cyclic} prompted our interest in the results of
this paper. The second-named author would like to thank
Universit{\'e} Catholique de Louvain for its hospitality during
the time this work was being conducted.

\section{Factorization of an Evaluation Fibration}\label{sec:w-factor}

The main purpose of this section is the proof of \thmref{thm:Gamma
factor}, \thmref{thm:Gamma+ factor} and their corollaries. The
results will flow from some general considerations about
fibrations of nilpotent spaces $F \to E \to X$ in which both $F$
and $E$ are $H_0$-spaces. The evaluation fibration $\omega\colon
\map(X,X;1) \to X$ is of this form with fibre the subspace of
$\map(X,X;1)$ consisting of based maps, which we denote by
$\map_*(X,X;1)$.

First we focus on the fibre inclusion of such a fibration.

\begin{proposition}\label{prop:H-space splitting}
Suppose $j \colon F \to E$ is any map between $H_0$-spaces.  Then
$E_\Q$ decomposes up to homotopy equivalence as $E_\Q \simeq Y \times
Z$, with $Y$ and $Z$ products of rational \EM spaces, so that there is
a corresponding map $\phi\colon Y \to F_\Q$ with
$(j_\Q)_\#\big(\pi_*(F_\Q)\big)= i_\#(\pi_*(Y))$ and $j_\Q\circ \phi =
i$, where $i\colon Y \to E_\Q$ denotes the inclusion of the first
factor.
\end{proposition}

\begin{proof}
Decompose $\pi_*(E_{\mathbb Q})$ as $V \oplus W$ with $V = \im
(j_{\Q})_\#$ and $W$ a complement.  Set $Y = \prod_i K(V_i,i)$ and
$Z = \prod_iK(W_i,i)$. Choose a basis $\{\alpha_t\}_{t \in T}$ for
$V$ so that $Y = \prod_t K(\Q, |\alpha_t|)$.  If $|\alpha_t|$ is
odd, then we identify $K(\Q,|\alpha_t|) \simeq
S^{|\alpha_t|}_{\Q}$ and we have a map $\alpha_t \colon
K(\Q,|\alpha_t|) \to E_\Q$.  If $|\alpha_t|$ is even, we construct
a corresponding map as follows.  First, we identify
$K(\Q,|\alpha_t|) \simeq \Omega \Sigma S^{|\alpha_t|}_{\Q}$.  Let
$\epsilon \colon X \to \Omega \Sigma X$ denote the adjoint of the
(suspension of the) identity.  Since $E_\Q$ is an $H$-space, we
may choose a retraction $r \colon \Omega \Sigma E_\Q \to E_\Q$ of
$\epsilon \colon E_\Q \to \Omega \Sigma E_\Q$ so that $r\circ
\epsilon = 1$ and the following diagram commutes:
$$\xymatrix{\Omega \Sigma S^{|\alpha_t|}_{\Q} \ar[r]^{\Omega \Sigma \alpha_t} & \Omega \Sigma
E_\Q \ar@/^1pc/[d]^{r}\\
S^{|\alpha_t|}_{\Q} \ar[u]^{\epsilon} \ar[r]_{\alpha_t} & E_\Q
\ar[u]^{\epsilon} }$$
That is, each $\alpha_t$ extends to map $\tilde{\alpha}_t = r\circ
\Omega \Sigma \alpha_t \colon K(\Q,|\alpha_t|) \to E_\Q$.  So far,
we have a map $a \colon \bigvee_t K(\Q,|\alpha_t|) \to E_\Q$
defined as $\alpha_t$ on odd-degree summands and
$\tilde{\alpha}_t$ on even-degree summands.  Now we may use the
multiplication of $E_\Q$ to extend this map to the product,
yielding a map $A \colon Y \to E_\Q$.  From the construction, we
have that $\im A_\# = V = \im (j_\Q)_\#$.

An identical construction yields a map $B \colon Z \to E_\Q$ that
satisfies $\im B_\# = W$.  Finally, one more use of the
multiplication $m$ of $E_\Q$ gives a map $m\circ(A\times B) \colon
Y \times Z \to E_{\mathbb Q}$ that is a homotopy equivalence.

For each $\alpha_t$, choose a $\beta_t \in \pi_*(F_\Q)$ such that
$(j_\Q)_\#(\beta_t) = \alpha_t$.  Repeating the above argument
with the $\beta_t$ replacing the $\alpha_t$ yields a map $\phi
\colon Y \to F_\Q$ with the desired properties, namely that
$j_\Q\circ\phi = A = m\circ(A\times B)\circ i$.
\end{proof}

Now consider any map $p\colon E \to X$ with $E$ an $H_0$-space. We
will construct a counterpart to the total Gottlieb element that
depends on the map $p$.  We first show that, under the hypothesis
that $X$ is finite---or more generally of finite rational
category, the image of $p_\#$ in rational homotopy groups is
restricted exactly as in the F{\'e}lix-Halperin result about
Gottlieb groups mentioned in the introduction.  Indeed, the
following result generalizes that result.  Here, we denote the
rational category of $X$ by $\cat0(X)$ (see \cite{F-H} or
\cite{F-H-T} for details of this invariant).

\begin{proposition}
Let $X$ be a nilpotent space and $p\colon E \to X$ be any map with
$E$ an $H_0$-space.  If $X$ has finite rational category, that is,
if $\cat0(X) = r < \infty$, then
$p_\#\big(\pi_{\mathrm{even}}(E_\Q)\big) = 0$ and
$p_\#\big(\pi_{\mathrm{odd}}(E_\Q)\big)$ is of (finite) dimension
no more than $r$.
\end{proposition}

\begin{proof}
For the first assertion, suppose that $\beta \in \pi_{2i}(E_\Q)$.
Because $E_\Q$ is an $H$-space the map $\beta$ extends, exactly as
in the proof of the previous result, to a map $\tilde\beta \colon
\Omega S^{2i+1} \to E_{\mathbb Q}$ that is injective in (rational)
homotopy in degree $2i$.  If $p_\#(\beta) \not= 0$, then
$p\circ\tilde\beta \colon \Omega S^{2i+1} \to X$ is a map that is
injective in rational homotopy---recall that $\Omega S^{2i+1}$
rationalizes to an \EM space $K(\Q, 2i)$. But then the mapping
theorem of \cite{F-H} implies that $\infty = \cat0(K(\Q, 2i)) \leq
\cat0(X) = r$, which is a contradiction. Therefore, we have
$p_\#\big(\pi_{\mathrm{even}}(E_\Q)\big) = 0$.  For the second
assertion, consider any finite, linearly independent subset
$\{\alpha_1, \ldots, \alpha_k\}$ of $\pi_{\mathrm{odd}}(X_\Q)$
such that each $\alpha_i \in \pi_{n_i}(X_\Q)$ is in the image of
$p_\#$.  Choose a $\beta_i \in \pi_{n_i}(E_\Q)$ with
$p_\#(\beta_i) = \alpha_i$ for each $i$. Write the corresponding
product of odd-dimensional rational spheres $\prod_{i=1}^k
S_\Q^{n_i}$ as $S_p$. Using the multiplication of $E_\Q$, we may
extend the map $\bigvee_i S_\Q^{n_i} \to E_\Q$ defined as
$\beta_i$ on each summand into a map $\widetilde{\Gamma}_p \colon
S_p \to E_{\mathbb Q}$. Specifically, let $M\colon E_\Q \times
\cdots \times E_\Q \to E_\Q$ denote the association $m\circ
(1\times m) \circ \cdots \circ (1\times \cdots \times 1 \times
m)$, where $m$ denotes the multiplication on $E_\Q$. (Recall that
we are not assuming $E_\Q$ to be associative.) Then we set
$$\widetilde{\Gamma}_p = M\circ (\beta_1\times\cdots\times\beta_k) \colon S_p \to E_{\mathbb
Q}.$$
Now an odd-dimensional rational sphere $S_\Q^{n_i}$ is a rational
\EM space $K(\Q ,{n_i})$.  From the construction, therefore, we
have that $p\circ \widetilde{\Gamma}_p \colon S_p \to X_{\mathbb
Q}$ is injective in rational homotopy groups.  Once again, the
mapping theorem implies that $k = \cat0(S_p) \leq r$.  The second
assertion follows.
\end{proof}

So now suppose that $p\colon E \to X$ is any map from an
$H_0$-space $E$ to a nilpotent, finite space $X$.  The image of
$p$ in rational homotopy groups is of finite dimension and we may
pick a finite basis $\{\alpha_1, \ldots, \alpha_k\}$ in
$\pi_{\mathrm{odd}}(X_\Q)$ for this image. Exactly as in the above
proof, we construct a map $\widetilde{\Gamma}_p \colon S_p \to
E_{\mathbb Q}$ and then set $\Gamma_p =
p_\Q\circ\widetilde{\Gamma}_p \colon S_p \to X_\Q$. (In the case
in which $p$ has trivial image in rational homotopy groups, we may
take $\widetilde{\Gamma}_p$ and $\Gamma_p$ to be the trivial map.)
In all cases, our construction gives a commutative diagram
$$
\xymatrix{ & E_\Q \ar[d]^-{p_\Q}\\ S_p
\ar[ru]^{\widetilde{\Gamma}_p} \ar[r]_{\Gamma_p}&  X_\Q}
$$
in which $\Gamma_p$ is both injective and onto the image of $p$ in
rational homotopy groups.

\begin{definition}\label{def:total Gottlieb}
Suppose given any map $p\colon E \to X$ from an $H_0$-space $E$ to
a nilpotent, finite space $X$.  A \emph{total Gottlieb element for
$X_\Q$ with respect to $p$} is a map $\Gamma_p \colon S_p \to
X_\Q$ that admits a lift $\widetilde{\Gamma}_p \colon S_p \to
E_\Q$ through $p_\Q$, where
\begin{enumerate}
\item $S_p$ is a product of rational \EM spaces with homotopy
isomorphic to $\im (p_\Q)_\#\colon \pi_*(E_\Q) \to \pi_*(X_\Q)$;
and
\item $\Gamma_p$ is injective in (rational) homotopy groups.
\end{enumerate}
\end{definition}

In general, there may be many choices of total Gottlieb elements
with respect to $p$ and different lifts of each. By the above
discussion, we see that such always exist.  We keep the notation
$\Gamma_X \colon S_X \to X_\Q$ for a total Gottlieb element with
respect to the ordinary evaluation fibration $\omega \colon
\map(X,X;1) \to X$.

\begin{theorem}\label{thm:Gamma factor general}
Let
$$\xymatrix{F \ar[r]^{j} & E \ar[r]^{p} & X}$$
be a fibration sequence of nilpotent spaces in which $F$ and $E$
are $H_0$-spaces and $X$ is a nilpotent, finite space. Let
$\Gamma_p \colon S_p \to X_\Q$ be any total Gottlieb element for
$X_\Q$ with respect to $p$ and $\widetilde{\Gamma}_p$ any lift of
of $\Gamma_p$ through $p_\Q$.  Assume there is an action $A \colon
F_{\Q} \times E_{\Q} \to E_{\Q}$ of $F_{\Q}$ on $E_{\Q}$ that
  satisfies $A\circ i_1 = j_\Q$ and $p_\Q\circ A = p_\Q\circ p_2
\colon F_\Q\times E_\Q \to X_\Q$. Then there is a retraction $r
\colon E_\Q \to S_p$ of $\widetilde{\Gamma}_p$ such that $p_\Q =
\Gamma_p\circ r\colon E_\Q \to X_\Q$.
\end{theorem}

\begin{proof}
From \propref{prop:H-space splitting}, we assume an identification
$E_\Q \simeq Y \times Z$, with $Y$ and $Z$ rational $H$-spaces,
together with maps $i\colon Y \to E_\Q$ and $\phi\colon Y \to
F_\Q$ with $i_\#$ an injection onto $\im (j_\Q)_\#$ and $j_\Q\circ
\phi = i$. Now consider the following commutative diagram:
$$\xymatrix{Y \times S_p \ar[r]^-{p_2}
\ar[d]_{A\circ(\phi\times\widetilde{\Gamma}_p)}^{\simeq} & S_p \ar[d]^{\Gamma_p}\\
E_\Q \ar@{.>}@/_1pc/[u]_{H}\ar[r]_-{p_\Q} & X_\Q}$$
Observe that $A\circ(\phi\times\widetilde{\Gamma}_p)\circ i_1 =
i\colon Y \to E_\Q$.  Furthermore, from the long exact sequence in
homotopy of the fibration, we find that
$A\circ(\phi\times\widetilde{\Gamma}_p)\circ i_2 \colon S_p \to
E_\Q$ has image in homotopy that is complementary to $\im
(j_\Q)_\#$.  Hence $A\circ(\phi\times\widetilde{\Gamma}_p)$
induces an isomorphism in rational homotopy and thus is a homotopy
equivalence. Consequently, there is an inverse (rational) homotopy
equivalence $H \colon E_\Q \to Y\times S_p$ as indicated in the
diagram.  Now set $r = p_2\circ H \colon E_\Q \to S_p$.  Then we
have $r \circ \widetilde{\Gamma}_p = p_2\circ H \circ
A\circ(\phi\times\widetilde{\Gamma}_p) \circ i_2 = 1 \colon S_p
\to S_p$, so that $r$ is a retraction of $\widetilde{\Gamma}_p$.
Furthermore, since $p\circ A(\varphi\times \widetilde{\Gamma}_p) =
\Gamma_p\circ p_2$, we have $\Gamma_p\circ r = \Gamma_p\circ
p_2\circ H = p_\Q \colon E_\Q \to X_\Q$, which gives the desired
factorization.
\end{proof}

We obtain \thmref{thm:Gamma factor} by specializing as follows:

\begin{proof}[Proof of \thmref{thm:Gamma factor}]
The action
$$A \colon \map_*(X,X;1)\times \map(X,X;1) \to
\map(X,X;1),$$
defined by $A(f,g) = g\circ f$, restricts to the inclusion
$\map_*(X,X;1)\to \map(X,X;1)$ and satisfies the hypothesis of
\thmref{thm:Gamma factor general}. Therefore, we may apply the
result to the evaluation fibration sequence
$$\map_*(X,X;1) \to \map(X,X;1)\stackrel{\omega}{\to} X$$
and the total Gottlieb element for this evaluation map constructed
from the Gottlieb groups as in the introduction.
\end{proof}

By the same argument, we obtain \thmref{thm:Gamma factor} for each of
the evaluation fibrations in which $\top$, $\diff$, and so-forth,
replaces $\map$, as in the introduction.

The following observation allows us to strengthen
\thmref{thm:Gamma factor} in certain circumstances.  We will also
use it in \secref{sec:cohomology}.  Roughly speaking, we may say
that if $X$ decomposes up to homotopy equivalence as a product,
then \emph{the evaluation map decomposes as a corresponding
product of evaluation maps}.

More precisely, suppose we have a homotopy equivalence $h\colon X
\to A \times B$. Then we have homotopy equivalences $$h_* \colon
\map(X, X;1) \to \map(X, A\times B;h)$$ and $$h^* \colon
\map(A\times B, A\times B;1) \to \map(X, A\times B;h)\,.$$ Let
$p_1\colon A \times B \to A$ and $p_2\colon A \times B \to B$
denote the projections, and $i_1\colon A \to A \times B$ and
$i_2\colon B \to A \times B$ the inclusions.   We write $h = (h_1,
h_2)$, with $h_1 = p_1\circ h$ and $h_2 = p_2\circ h$, and define
evaluation maps $\omega_1 \colon \map(X,A; h_1) \to A$ and
$\omega_2 \colon \map(X,B; h_2) \to B$ by evaluation at the
basepoint of $X$.  Let $I \colon \map(A\times B, A\times B; 1) \to
\map(A\times B, A; p_1) \times \map(A\times B, B; p_2)$ be the
standard homeomorphism. One checks easily that
$$h\circ \omega_X = (\omega_1 \times
\omega_2)\circ I\circ (h^*)^{-1} \circ h_* \colon \map(X, X;1) \to
A\times B,$$
and thus we may identify $\omega_X$ with $h^{-1}\circ(\omega_1
\times \omega_2)\circ I \circ h_*$.  Furthermore, we observe that
$\omega_1$ factors through $\omega_A$ as $\omega_1 = \omega_A
\circ (i_1)^*$ and likewise $\omega_2 = \omega_B \circ (i_2)^*$.
Therefore, we have the following commutative diagram:
\begin{equation}\label{eq:w product}
\xymatrix{ \map(X, X;1) \ar[dd]_-{\omega_X} \ar[rr]^-{I\circ
(h^*)^{-1}\circ h_*}_-{\simeq}  & & \map(A\times B, A; p_1) \times
\map(A\times B, B; p_2) \ar[d]^-{(i_1)^* \times
(i_2)^*}\\
   & & \map(A, A; 1_A)
\times \map(B, B; 1_B) \ar[d]^-{\omega_A \times \omega_B}\\
X \ar[rr]^-{\simeq}_-{h} & & A\times B.}
\end{equation}

This discussion leads to the following result, which should be
compared with the well-known fact that $G_*(A\times B) \cong
G_*(A) \oplus G_*(B)$ \cite{Go1}.

\begin{theorem}\label{thm: X = A x B}
Suppose that we have a homotopy equivalence $X \simeq A \times B$.
Then the evaluation map $\omega_X$ factors through the product of
evaluation maps $\omega_A \times \omega_B$.
\end{theorem}

We now continue with the main results. In order to study
generalized evaluation maps $\ev \colon E \to X$, we first present
a global structure result concerning maps between $H_0$-spaces.

\begin{proposition}\label{prop:model of H-map}
Let $f \colon X \to Y$ be a map between $H_0$-spaces.
\begin{enumerate}
\item[(a)] The map $f$ admits a Sullivan
minimal model of the form $\varphi \colon (\land (V\oplus R),0)
\to (\land(V\oplus S),0)$ with $\varphi(v) = v$ for $v \in V$ and
such that $\varphi(R) \in \land^{\geq 2}(V \oplus S) \cap \land V
\otimes \land^+(S)$.
\item[(b)] If $f_\Q$ is an $H$-map then $f$ admits a model of the
form $\varphi \colon (\land(V\oplus K),0) \to (\land(V\oplus
S),0)$ with $\varphi(v) = v$ for $v \in V$ and $\varphi(K) = 0$.
\end{enumerate}
\end{proposition}

\begin{proof} Let $\varphi \colon (\land T,0) \to (\land W,0)$ be
any model of $f$.  We will use standard tricks from rational
homotopy to change generators in $\land T$ and $\land W$ so that,
with respect to the new generators, the minimal model of $f$ has
the desired form.

(a) We denote by $V$ a maximal subspace of $T$ such that
$Q(\varphi) \colon V \to W$ is injective. Denote by $R\subseteq T$
a complement of $V$ and by  $S \subseteq W$ a complement of $\im
Q(\varphi)$ in $W$. Let $\{v_i\}_{i\in I}$ be a graded basis for
$V$. Then the elements $\varphi (v_i)$ are linearly independent
indecomposable elements in $\land W$.  Denote by $\{r_j\}_{j\in
J}$ a graded basis for $R$ and $\{s_k\}_{k\in K}$ a graded basis
for $S$. With respect to the generators $\{v_i, r_j\}$ for $\land
T$ and $\{v'_i = \varphi (v_i), s_k\}$ for $\land W$, the map
$\varphi$ satisfies $\varphi (v_i)=v'_i$ and $\varphi (R) \subset
\land^{\geq 2}(W)$. We can thus suppose $\varphi(v) = v$ and that
$\varphi (R)$ is decomposable. We now change generators in $R$ so
that $\varphi (R)$ also belongs to the ideal generated by $S$.
Suppose that this is true for $R^{<n}$, and let $r$ be a generator
in $R^n$. If $\varphi (r) = a + b$ with $a\in \land V$ and $b$ in
the ideal generated by $S$, we change the generator to $r' = r -
a$.  The result follows by induction.

(b) Here, we apply the previous step to write $\varphi \colon
\land(V \oplus K) \to \land(V\oplus S)$ with $\varphi(v) = v$ for
$v\in V$ and $\varphi(k)$ both decomposable and in $\land V
\otimes \land^+(S)$ for $k \in K$. We now prove by induction that
$\varphi$ is zero on $K$.

The existence of multiplications on $X_\Q$ and $Y_\Q$ is reflected
in their Sullivan models by morphisms of algebras $\Delta_1 \colon
\land T \to \land T \otimes \land T$ and $\Delta_2 \colon \land W
\to \land W \otimes \land W$ that satisfy $\Delta_1(v) -(v\otimes
1 + 1\otimes v) \in \land^+T\otimes \land^+T$ and likewise for
$\Delta_2$. Furthermore, since $f_\Q$ is an $H$-map, we have the
following commutative diagram after the previous step:
$$\xymatrix{ \land(V\oplus K) \ar[r]^-{\Delta_1} \ar[d]_{\varphi} &
\land(V\oplus K)\otimes\land(V\oplus K) \ar[d]^{\varphi\otimes
\varphi}\\ \land(V\oplus S) \ar[r]_-{\Delta_2} & \land(V\oplus S)
\otimes\land(V\oplus S)} $$
Assume inductively that we have $\varphi(K^{\leq n})=0$ and
suppose that $u\in K^{n+1}$.  We write
$$\varphi(u) = \varphi_r(u) +\varphi_{r+1}(u) + \cdots +
\varphi_m(u)$$
with $\varphi_r(u) \in \land^r(V\oplus S)$. By the definition of
$K$, we have $r\geq 2$.  Consider a term in $\varphi_r(u)$ that is
of minimal length $q$ in $\land S$, for some $1 \leq q \leq r$.
Let $\{s_i\}$ be a basis of $S$ and write such a minimal term as
$s_{i_1}s_{i_2}\cdots s_{i_q}\nu$ for some $\nu \in \land^{r-q}V$.
Then $\Delta_2 \varphi(u)$ contains a contribution $s_{i_1}\otimes
s_{i_2}\cdots s_{i_q}\nu$ and this term will appear uniquely as
such in $\Delta_2\varphi(u) - (1\otimes \varphi(u) +
\varphi(u)\otimes1)$. On the other hand, $\Delta_1(u) - (1\otimes
u + u\otimes1) \in \land^{+}(V \oplus K^{\leq n})\otimes
\land^{+}(V \oplus K^{\leq n})$ and so $(\varphi\otimes\varphi)
\Delta_1(u) - (1\otimes\varphi(u) + \varphi(u)\otimes1)$ cannot
contain any occurrence of a term such as $s_{i_1}\otimes
s_{i_2}\cdots s_{i_q}\nu$, by our induction hypothesis. In
summary, if $\varphi_r(u)$ contains some non-zero term, then we
cannot have $(\varphi\otimes\varphi)\Delta_1(u) = \Delta_2
\varphi(u)$, which is a contradiction.  It follows by induction
that $\varphi(K)=0$.
\end{proof}

We remark in passing that \propref{prop:model of H-map} implies the
following result:

\begin{corollary}
Let $f \colon X \to Y$ be a map between $H_0$-spaces that is an $H$-map
after rationalization. If $(f_\Q)_\#$ is zero, then $f$ is rationally
null-homotopic.
\end{corollary}

We also observe that the conclusion of \propref{prop:model of
H-map} (b) holds for certain compositions.  We will use this
observation in the following form in the sequel,:

\begin{corollary}\label{cor:composition of H-map and surj}
Let $g \colon X \to Y$ and $r \colon Y \to Z$ be maps between
$H_0$-spaces. If $g_\Q$ is an $H$-map and $(r_\Q)_\#$ is
surjective, then their composition $r\circ g$ admits a Sullivan
minimal model of the form $\varphi \colon (\land(V\oplus K), 0)
\to (\land(V\oplus W),0)$ with $\varphi(v) = v$ for $v \in V$ and
$\varphi(K) = 0$.
\end{corollary}

\begin{proof}
Denote by $\varphi \colon (\land W_1,0) \to (\land W_2,0)$ a
minimal model of $r\circ g$, and by $V \subset W_1$ a maximal
subspace such that $Q(\varphi) \colon V \to W_2$ is injective.
Then by part (a) of \propref{prop:model of H-map}, we have models
for $r$ and $g$
$$(\land (V\oplus R,0) \stackrel{\theta_1}{\longrightarrow} (\land
(V \oplus R \oplus S), 0) \stackrel{\theta_2}{\longrightarrow}
(\land (V \oplus W),0)\,,$$
with $\theta_1(v) = \theta_2(v) = v$. Using part (b) of
\propref{prop:model of H-map}, we can suppose that $\theta_2(R
\oplus S) = 0$. By a change of generators in $(\land (V\oplus
R),0)$ we can suppose that $\theta_1(R)$ is contained in the ideal
generated by $R \oplus S$, so that $\theta_2\circ\theta_1(R) = 0$.
\end{proof}

We now proceed to the proof of  our second main result, namely
\thmref{thm:Gamma+ factor},

\begin{proof}[Proof of \thmref{thm:Gamma+ factor}]
Suppose $\ev\colon E \to X$ is an evaluation map.  Then there is an
action $A \colon E \times X \to X$ that restricts to $\ev$.  The
adjoint $g\colon E \to \map(X,X;1)$ of this action, defined by
$g(y)(x) = A(y,x)$, is a lift of $\ev$ through $\omega \colon
\map(X,X;1) \to X$.  Since we assume the action is associative, the
adjoint $g$ is an $H$-map.  Upon rationalizing, we obtain the
commutative diagram
$$\xymatrix{E\ar[r]^-{g} \ar[dr]_{\ev_\Q} & \map(X,X_\Q; e)
\ar[r]^-{r}\ar[d]^{\omega_\Q} & S_X \ar[dl]^{\Gamma_X}\\
 &   X_\Q }$$
in which $r \colon \map(X,X_\Q; e) \to S_X$ is a retraction of
$\widetilde{\Gamma}_X \colon S_X \to \map(X,X_\Q; e)$ as in
\thmref{thm:Gamma factor}.  Since $r$ is a retraction, $r_\#$ is
surjective.  So we may apply \corref{cor:composition of H-map and
surj} and assume a model of $r\circ g \colon E \to S_X$ has the
form $\varphi \colon (\land(V\oplus K), 0) \to (\land(V\oplus
W),0)$ with $\varphi(v) = v$ for $v \in V$ and $\varphi(K) = 0$.
Thus $\varphi$ factors in the form
$$\xymatrix{\land(V\oplus K) \ar[r]_-{\mathrm{proj}} & \land V \ar[r]_-{\mathrm{incl}} & \land(V\oplus W) \ar@/_1pc/[l]_-{\mathrm{proj}}}$$
together with the evident retraction of the inclusion $\land V \to
\land(V\oplus W)$ as indicated.  When translated into spaces, this
implies that $r\circ g$ factors rationally through a rational
$H$-space $Y$
$$\xymatrix{E \ar[r]_{q} & Y \ar[r]^{j} \ar@/_1pc/[l]_{i} & S_X}.$$
Notice that $j \colon Y \to S_X$ has minimal model the projection
$\land(V\oplus K) \to \land V$ and hence is injective in
(rational) homotopy.  Furthermore, we have the right inverse $i$
for $q$ as indicated.   That is, we have maps that satisfy $j\circ
q = r \circ g$ and $q\circ i = 1$. Now consider the diagram
$$
\xymatrix{ & E \ar[d]^-{\ev_\Q} \ar@/_1pc/[ld]_{q}\\ Y \ar[ru]_{i}
\ar[r]_{\Gamma_X\circ j}&  X_\Q,}
$$
in which we have $\Gamma_X\circ j \circ q = \Gamma_X\circ r \circ
g = \omega_\Q \circ g = \ev_\Q$ and hence $\ev_\Q\circ i =
\Gamma_X\circ j \circ q\circ i = \Gamma_X\circ j$.  We see that
$\Gamma_X\circ j \colon Y \to X_\Q$ satisfies the requirements of
a total Gottlieb element for $X_\Q$ with respect to $\ev$.  Since
we have a retraction $q$ of $i$, which here serves as our lift of
$\Gamma_X\circ j$ through $\ev_\Q$, this total Gottlieb element
satisfies the conclusion of the theorem.

The conclusion now follows for every total Gottlieb element For
suppose given another total Gottlieb element    $\Gamma'_p\colon
S'_p \to X_\Q$ for $X_\Q$ with lift   $\widetilde{\Gamma}'_p
\colon S'_p \to E$.    Then the map $h =
q\circ\widetilde{\Gamma}'_p \colon S'_p \to Y$ is a homotopy
equivalence.  This follows since $S'_p$ and $Y$ have isomorphic
(rational) homotopy groups, and $h$ is injective in (rational)
homotopy groups. Therefore, we may define $r' = h^{-1}\circ
q\colon E \to S'_p$, which is easily checked to be a retraction of
$\widetilde{\Gamma}'_p$ that satisfies $\Gamma'_p\circ r' =
\ev_\Q$.
\end{proof}

We may supplement the vocabulary of \defref{def:total Gottlieb} with
the following:  Suppose given any map $p\colon E \to X$ from an
$H_0$-space $E$ to a nilpotent, finite space $X$. Then we define the
\emph{$n$th Gottlieb group of $X$ with respect to $p$} as the
subgroup of $\pi_n(X)$ that is the image of $p_\# \colon \pi_n(E)
\to \pi_n(X)$.  We denote this subgroup by $G^p_n(X)$.  Then we have
\corref{cor:w-sharp0}, phrased using this notation, as an immediate
consequence of \thmref{thm:Gamma+ factor}.

\begin{corollary}[\corref{cor:w-sharp0}]\label{cor:G ev zero}
Let $X$ be any nilpotent, finite complex and $\ev \colon E \to X$ an
evaluation map. Then $G_*^\ev(X_\Q) = 0$ if and only if $\ev_\Q$ is
null-homotopic.
\end{corollary}

Before we present some examples, we notice the following
generalization of \corref{cor:w-sharp0} that  does not require the
fibration to be ``principal" in the sense required by
\thmref{thm:Gamma factor general}:

\begin{theorem}\label{thm:w-sharp0}
Suppose given any fibration sequence of nilpotent spaces
$$\xymatrix{F \ar[r]^{j} & E \ar[r]^{p} & B}$$
in which $F$ and $E$ are $H_0$-spaces.  If $(p_\Q)_\# = 0\colon
\pi_*(E_\Q) \to \pi_*(B_\Q)$, then $p_\Q = *$.
\end{theorem}

\begin{proof}
From the long exact sequence in rational homotopy groups induced
by the fibration sequence, we have that $(j_\Q)_\#\colon
\pi_*(F_\Q) \to \pi_*(E_\Q)$ is surjective.    This gives  a
section $\sigma  \colon E_\Q \to F_\Q$ of the rationalized fibre
inclusion $j_\Q \colon F_\Q \to E_\Q$. Thus we have $p_\Q =
p_\Q\circ j_\Q\circ\sigma = *$, since $p_\Q\circ j_\Q = *$.
\end{proof}

\begin{example}
We give first an example of a fibration that satisfies the
hypotheses of \thmref{thm:Gamma factor general}, and yet is not a
cyclic map and therefore, in particular, does not satisfy the
hypotheses of \thmref{thm:Gamma+ factor}.  For this, let $B$
denote a space whose minimal model is $\Lambda(a, b, c)$, with
$|a| = |b| = 3$ and $|c| = 5$, and differential given by $d(a) =
d(b) = 0$ and $d(c) = ab$.  Then consider the map $p\colon S^3 \to
B$ that corresponds to one of the homotopy elements of $\pi_3(B)$.
We find that, up to rational equivalence, the homotopy fibre of
$p$ is the $H$-space $F = \Omega(S^3\times S^5)$.  Furthermore,
again up to rational equivalence, the fibre inclusion $j\colon F
\to S^3$ is null-homotopic.  The fibre sequence $F \to S^3 \to B$,
therefore, admits an action of $F$ on $S^3$ that is principal in
the sense required by the hypotheses of \thmref{thm:Gamma factor
general}. Namely, the projection $p_2\colon F\times S^3 \to S^3$
is such an action.  Observe, however, that the fibre map $p\colon
S^3 \to B$ cannot be a cyclic map, since $G_3(B) = 0$.  In
particular, this example does not satisfy the hypotheses of
\thmref{thm:Gamma+ factor}.
\end{example}

\begin{example}\label{ex:not partial}
Observe, however, that there are maps from an $H_0$-space that
induce zero on (rational) homotopy groups, and yet are not
(rationally) null-homotopic.  For instance, the quotient map $q
\colon S^3 \times S^3 \to S^6$ is a map from an $H$-space that
induces zero on rational homotopy groups, yet is non-zero on
rational cohomology groups and so is not rationally trivial. Of
course, here the homotopy fibre of $q$ is not an $H_0$-space.  It
is interesting to note that \corref{cor:G ev zero} implies $q_\Q$
\emph{cannot occur as the connecting map of any fibration}
(cf.~\corref{cor:partial trivial}).
\end{example}

Allowing a non-trivial image in homotopy for $p$ appears to make a
fundamental change in the situation.  In particular, if we simply
assume $F$ and $E$ are $H_0$-spaces, as in \thmref{thm:w-sharp0},
but allow the image of $p$ in rational homotopy groups to have
dimension $1$, then it may be impossible to factor $p$ through an
odd-dimensional sphere, or any finite product of odd-dimensional
spheres. We give examples to illustrate this point:

\begin{example}
Let $q \colon S^3\times S^3 \times S^3 \to S^9$ be the map
obtained by pinching out all but the top cell of the product. As
 may be checked by a direct computation, the fibre sequence
$$\xymatrix{F \ar[r]^-{j} & S^3\times S^3 \times S^3 \ar[r]^-{p} &
S^3\times S^9}$$
with $p = (p_1, q)$ has fibre that is rationally equivalent to the
$H$-space $S^3 \times S^3 \times K(\Q, 8)$.  Hence, the fibre
inclusion $j$ is a map of $H_0$-spaces. Now $p$ has image of
dimension $1$ on rational homotopy groups. Evidently, however, $p$
does not factor through $S^3$ (or any single odd-dimensional
sphere).
\end{example}

\begin{example}
We describe a rational fibre sequence of  the form
$$\xymatrix{F \ar[r]^-{j} & E \ar[r]^-{p} & S^3\vee S^9},$$
in which $E$ and $F$ are $H_0$-spaces and where $p$ does not
factor through any finite product of odd-dimensional spheres.
First we specify a map of minimal models $\mathcal{M}_{p} \colon
\mathcal{M}_{S^3\vee S^9} \to \mathcal{M}_{E}$ by writing
$\mathcal{M}_{S^3\vee S^9} = \Lambda(b, \{u_i\}_{i\geq 1};d)$,
with $|b|= 3$ and $|u_i| \geq 9$, and then setting $\mathcal{M}_E$
to be the minimal model $\Lambda(b, y, \{v_i\}_{i\geq 1}; d_E=0)$,
with $|b| = |y| = 3$ and $|v_i| = |u_i| - 6$. The map
$\mathcal{M}_p$ is given by $\mathcal{M}_p(b) = b$, and
$\mathcal{M}_p(u_i) = byv_i$ for $i \geq 1$.  Since $d$ is
decomposable, we have $\mathcal{M}_p\circ d = 0$ and thus
$\mathcal{M}_p$ is a map of DG algebras.  Hence it defines a map
of rational spaces $p \colon E \to (S^3\vee S^9)_\Q$.  By standard
rational homotopy techniques, one checks that the homotopy fibre
of $p$ is an $H_0$-space. However, one may see from the minimal
models that the map $p$ does not factor through any finite product
of odd-dimensional spheres.
\end{example}

\section{Gottlieb groups and homotopy monomorphisms}\label{sec:technical}

Let $\ev \colon E \to X$ be an evaluation map.  By
\thmref{thm:Gamma+ factor}, $\ev$ factors as $\ev = \Gamma_w\circ
r$ where $r \colon E_\Q \to S_{\ev}$ is a left inverse of
$\widetilde{\Gamma_\ev}$.   As a retraction, $r$ has
$\widetilde{\Gamma}_\ev$ as a right inverse and so \emph{is a
homotopy epimorphism}.  That is, the map of homotopy sets
$$r_* \colon [A, E_\Q] \to [A, S_\ev]$$
is surjective for any space $A$.  On the other hand, a total
Gottlieb element $\Gamma_\ev \colon S_\ev \to X_\Q$ generally does
not admit a left inverse.  For instance, take $X = S^2$ so that
$G_*(X_\Q) = G_3(S^2_\Q) \cong \Q$.  Then $S_X = S^3_\Q$ and we may
take $\Gamma_X \colon S_X \to X_\Q$ to be the rationalized Hopf map,
which does not admit a left inverse.  Nonetheless, we will show that
$$(\Gamma_\ev)_* \colon [A, S_\ev] \to [A, X_\Q]$$
is injective for any nilpotent space $A$
. In order to show this, and
in addition to obtain our results about cohomology, we need to
establish some technical points concerning Gottlieb groups and
rational homotopy monomorphisms.

The following discussion will fix our notation for the remainder of
the paper. Suppose $X$ has minimal model $(\land W, d_X)$. The
Gottlieb group $G_*(X_\Q)$ may be identified with the subspace of
$\Hom(W, \Q)$ formed by those linear maps that extend to derivations
of $\land W$ that commute with $d_X$ (see \cite{F-H-T} for a
discussion of this). Denote by $\overline{\theta}_i$ a linear basis
of $G_*(X_\Q)$, and by $v_i$ elements of $W$ with
$\overline{\theta}_i(v_j) = \delta_{ij}$.  We denote by $\theta_i$
an extension of $\overline{\theta}_i$ to a derivation of $\land W$
that satisfies $d_X \theta_i = (-1)^{|v_i|}\theta_i d_X  $. We
suppose, without loss of generality, that $|v_i| \leq |v_j|$ for $i
< j$.  Then we may---and do---suppose that $\theta_i(v_j) = 0$ for
$i > j$. Other than this, however, we have very little control over
how the $\overline{\theta}_i$ extend.  This point is the main source
of the technicalities.  We denote by $V$ the vector space generated
by the $v_i$, and $Z$ a choice of complement in $W$. Thus the
minimal model of $X$ is $(\land(V \oplus Z), d_X)$ with $V = \langle
v_1, \dots, v_r \rangle$ corresponding to the Gottlieb group,
accompanying derivations $\theta_1, \dots, \theta_r$, and $Z$ a
complement to $V$ in $W$.

\begin{lemma}\label{lem:Z theta stable}
With notation as above, we may choose $Z$, $V$, and the $\theta_i$
such that $\theta_i (Z \oplus \langle v_{i+1}, \cdots , v_r\rangle)
\subset \land V \otimes \land^+ Z$. In particular the ideal $I(Z)$
is $\theta_i$-stable for each $i$.
\end{lemma}

\begin{proof}
Let $\mathcal{L}$ denote the Lie algebra of derivations of
$\mathcal{M}_X$ generated by the derivations $\theta_1, \dots,
\theta_r$. We prove by induction on $k$ that we may choose
  $Z$ and $V$ for which we have $\theta(W) \subseteq \mathbb Q \oplus  (\land^{\geq
k}V + \land V\otimes \land^{+}Z)$ for any $\theta \in
\mathcal{L}$, for all $k$.  Since $\land V$ is finite dimensional,
taking $k > r$ establishes the result.

For $k=1$, we choose $Z = \cap_{i=1}^r \ker(\overline{\theta}_i
\colon W \to \Q)$.  We have directly that $\theta_i(Z) \subseteq
\land^+(V\oplus Z)$, and hence $\theta(Z) \subseteq
\land^+(V\oplus Z)$ for any $\theta \in \mathcal{L}$.

Now suppose that, for some $k\geq 1$, we have $\theta(W) \subseteq
\mathbb Q \oplus (\land^{k}V + \land V\otimes \land^+ Z)$ for any
$\theta \in \mathcal{L}$.  For each generating derivation
$\theta_j$, and for $z \in W$ a basis element, with $\theta (z)
\not\in \mathbb Q$,
 we write
$$\theta_j(z) \equiv \sum_{i_1<i_2<\cdots<i_k}
\;\lambda^{(i_1,i_2,\dots,i_k)}_j\;v_{i_1}v_{i_2}\cdots v_{i_k}$$
modulo terms in $\land^{\geq k+1}V + \land V\otimes \land^+ Z$.
Then we make a change of basis for $W$---in effect, a different
choice of complement---by replacing each basis element $z$ with
$z'$, where
$$z' = z -  \sum_{j=k+1}^r\,\sum_{i_1<i_2<\cdots<i_k<j}
\;\lambda^{(i_1,i_2,\dots,i_k)}_j\;v_j v_{i_1}v_{i_2}\cdots
v_{i_k}.$$
The effect of this basis change in   $W$ is that we may now
suppose
\begin{equation}\label{eq:changed Z}
\theta_j(z) \equiv \sum_{i_1<i_2<\cdots<i_k| i_k\geq j}
\;\lambda^{(i_1,i_2,\dots,i_k)}_j\;v_{i_1}v_{i_2}\cdots v_{i_k}
\end{equation}
modulo terms in $\land^{\geq k+1}V + \land V\otimes \land^+ Z$,
for each generating derivation $\theta_j$ and each element $z \in
W$ such that $\theta (z)\not\in \mathbb Q$. We now claim that all
the coefficients $\lambda^{(i_1,i_2,\dots,i_k)}_j$ that appear in
(\ref{eq:changed Z}) are in fact zero.  For suppose that this is
not the case, and let $j$ be the least index for which some
$\lambda^{(i_1,i_2,\dots,i_k)}_j$ in (\ref{eq:changed Z}) is
non-zero.  Denote by $n\geq j$ the maximum of the $i_k$ with
$\lambda^{(i_1,i_2,\dots,i_k)}_j \not=0$.  Then
$\theta_n\circ\theta_j(z) = \alpha + \beta$, with $\alpha \not=0
\in \land^{r-1}(v_{i_1}, v_{i_2}, \dots, v_{n-1})$ and $\beta \in
\land^{\geq r}V + \land V\otimes \land^+ Z$.  If $n = j$, then
$\theta_n\circ\theta_j = \frac{1}{2} [\theta_n, \theta_j] \in
\mathcal{L}$, and this contradicts the induction hypothesis on
$\mathcal{L}$.   However, if $n > j$, then
$\theta_j\circ\theta_n(z) = \gamma + \delta$, with $\gamma$ of
length $k-1$ but in $\land(v_{i_1}, v_{i_2}, \dots, \widehat{v}_j,
\dots, v_{n-1})\otimes \land^+(v_{n}, v_{n+1}, \dots, v_{r})$ and
$\delta \in \land^{\geq r}V + \land V\otimes \land^+ Z$. This
shows again that $[\theta_n, \theta_j] \in \mathcal{L}$
contradicts the induction hypothesis on $\mathcal{L}$.  It follows
that all the coefficients $\lambda^{(i_1,i_2,\dots,i_k)}_j$ that
appear in (\ref{eq:changed Z}) are zero.  Therefore, we have
$\theta_j(W) \subseteq \mathbb Q \oplus ( \land^{k+1}V + \land
V\otimes \land^+ Z)$ for any $z \in W$, for each generating
derivation $\theta_j$.

To complete the inductive step, we must also consider a general
$\theta \in \mathcal{L}$.  Suppose that, for some $z\in W$, we
have
$$\theta(z) \equiv \sum_{i_1<i_2<\cdots<i_k}
\;\mu^{(i_1,i_2,\dots,i_r)}\;v_{i_1}v_{i_2}\cdots v_{i_k}$$
modulo terms in $\land^{\geq k+1}V + \land V\otimes \land^+ Z$. We
claim that all the coefficients $\mu^{(i_1,i_2,\dots,i_k)}$ that
appear in this expression are zero.  For suppose not, and once
again, denote by $n$ the maximum of the $i_k$ for which some
$\mu^{(i_1,i_2,\dots,i_k)} \not= 0$. The composition
$\theta_n\circ\theta(z)$ then contains a non-zero term in
$\land^{r-1}V$.  On the other hand, since $\theta$ is a
derivation, and we have just shown that $\theta_n(z) \in
\land^{\geq k+1}V + \land V\otimes \land^+ Z$, we have
$\theta\circ\theta_n(z) \in \land^{\geq k+1}V + \land V\otimes
\land^+ Z$.  Therefore, $[\theta_n, \theta] \in \mathcal{L}$
contradicts the induction hypothesis on $\mathcal{L}$. This shows
that all the $\mu^{(i_1,i_2,\dots,i_k)} = 0$, and hence the
induction is complete.
\end{proof}

\begin{proposition}\label{prop:I(Z) d-stable}
Suppose $V$, $Z$, and the $\theta_i$ satisfy $\theta_i(W) \subseteq
\mathbb Q \oplus (\land V\otimes \land^{+}Z)$ for each $i$. Then,
\begin{enumerate}
\item
$d_X(W) \subseteq \land V \otimes \land^{\geq 2}Z$.  In
particular, the ideal $I(Z)$ itself is $d_X$-stable.
\item There exists a choice of total Gottlieb element $\Gamma_X
\colon S_X \to X_\Q$ with minimal model $\mathcal{M}_{\Gamma} \colon
(\land(V \oplus Z), d_X) \to (\land V, d=0)$ that satisfies
$\mathcal{M}_{\Gamma}(Z) = 0$ and $\mathcal{M}_{\Gamma}(v) = v$ for
$v \in V$
.
\end{enumerate}
\end{proposition}

\begin{proof}
(1) We first argue by contradiction to prove that $d_X(W) \subset
\land V \otimes \land^+Z$. Suppose this is not true, and that $m
\geq 1$ is the minimal length for which any $d(\chi)$ contains a
non-zero term in $\land^{m}V$.  For such a $\chi \in \land
V\otimes \land Z$, write $d(\chi) = \alpha + \beta$ with $\alpha
\not= 0 \in \land^{m}V$ and $\beta \in \land^{\geq m+1}V + \land V
\otimes \land^{+}Z$. Further, suppose that $\alpha \in
\land^{m}(v_1, \dots, v_s)$ for some $s\leq r$ such that
$$d(\chi) = \alpha' + \alpha'' v_s + \beta,$$
with $\alpha' \in \land^{m}(v_1, \dots, v_{s-1})$ and
$\alpha''\not= 0 \in \land^{m-1}(v_1, \dots, v_{s-1})$ ($\alpha''
\in \Q$ if $m=1$). Then $\theta_sd(\chi) = \pm \alpha'' +
\theta_s(\beta)$ (recall that $\theta_i(v_j) = 0$ for $i > j$).
However, we have $\theta_sd(\chi) = - \theta_s d(\chi)$.  Using
\lemref{lem:Z theta stable} and the fact that $\theta_s$ is a
derivation, we also have $\theta_s(\beta) \in \land^{\geq m}V +
\land V \otimes \land^{+}Z$.  This contradicts our minimal length
assumption.

We claim that $d_X(W) \subset \land V \otimes \land^{\geq 2} Z$.
Suppose this is not the case and let $w$ be an element of lowest
degree such that $$d_X(w) = \sum_{i=1}^q z_i\omega_i + \alpha\,,$$
with $z_i \in Z$, $\vert z_1\vert \leq \vert z_2\vert \leq \cdots
\leq \vert z_q\vert$, $\omega_i \in \land V$ and $\alpha \in
\land^{\geq 2}Z \otimes \land V$.  We choose then an element $v_s$
of highest degree such that $\omega_q = v_s\gamma + \delta$,
$\gamma \neq 0$, $\gamma, \delta \in \land (v_1, \cdots ,
v_{s-1})$. Then $$\theta_qd_X(w) = z_qv_s \, \mbox{mod}\, \land V
\otimes (\land^{\geq 2} Z + Z^{<\vert z_g\vert} + (z_1, \cdots ,
z_{q-1}))\,.$$ Since $\theta_qd_X(w) = d_X \theta_q(w)$, there
exists an element $w'\in W$ with $\vert w'\vert <\vert w\vert $
such that $d_X(w') \not\in \land V \otimes ^{\geq 2} Z$. This is
impossible by our assumption.

 (2)
We will define a map $\phi \colon \mathcal{M}_{X} \to
\mathcal{M}_{S_X}\otimes\mathcal{M}_{X}$ whose composition with
the projection onto the first factor
$$(1\cdot\epsilon)\circ\phi\colon \mathcal{M}_X \to \mathcal{M}_{S_X}\otimes\mathcal{M}_X
\to \mathcal{M}_{S_X}$$
is surjective and satisfies $(1\cdot\epsilon)\circ\phi(Z) = 0$,
and whose composition with the projection onto the second factor
is the identity, $(\epsilon\cdot1)\circ\phi = 1\colon
\mathcal{M}_X \to \mathcal{M}_X$.

Translating this into topological terms, $\phi$ is the minimal
model of a map $F\colon S_X \times X_\Q \to X_\Q$ such that
$F\circ i_1 \colon S_X \to X_\Q$ is injective in rational homotopy
and $F\circ i_2 = 1 \colon X_\Q \to X_\Q$.  In other words, we may
choose $F\circ i_1$ as a total Gottlieb element (the corresponding
lift through $\omega_\Q$ is given by the adjoint of $F$).
Furthermore, the model of $F\circ i_1$ is
$(1\cdot\epsilon)\circ\phi$ by construction, which satisfies
$(1\cdot\epsilon)\circ\phi(Z) = 0$.

So as to avoid confusion, we write
$\mathcal{M}_{S_X}\otimes\mathcal{M}_{X}$ as $\land V'\otimes
\land V\otimes\land Z$, with $V' = \langle v'_1, \dots, v'_r
\rangle$. First, define a sequence of maps $\phi_1, \dots, \phi_r
\colon \mathcal{M}_{X} \to
\mathcal{M}_{S_X}\otimes\mathcal{M}_{X}$ by $\phi_1(\chi) = \chi +
v'_1\theta_1(\chi)$, and
$$\phi_s(\chi) = \phi_{s-1}(\chi) +
v'_s\theta_s\big(\phi_{s-1}(\chi)\big)$$
for $s = 2, \dots, r$.  Then we set $\phi = \phi_r$.  An inductive
argument shows that $\phi$ so defined is a DG algebra map.  For it
is straightforward to check that $\phi_1$ is a DG algebra map.
Supposing inductively that $\phi_{s-1}$ is a DG algebra map, the
computation
$$\begin{aligned}
\phi_{s-1}(\chi_1)\phi_{s-1}(\chi_2) & = (\phi_{s-1}(\chi_1) +
v'_s\theta_s\big(\phi_{s-1}(\chi_1)\big))(\phi_{s-1}(\chi_2) +
v'_s\theta_s\big(\phi_{s-1}(\chi_2)\big))\\
&= \phi_{s-1}(\chi_1)\phi_{s-1}(\chi_2) +
v'_s\theta_s\big(\phi_{s-1}(\chi_1)\big)\phi_{s-1}(\chi_2)\\
& \hskip1.5truein    + (-1)^{|\chi_1|}v'_s
\phi_{s-1}(\chi_1)\theta_s\big(\phi_{s-1}(\chi_2)\big)\\
&= \phi_{s-1}(\chi_1\chi_2) +
v'_s\theta_s\big(\phi_{s-1}(\chi_1\chi_2)\big)
\end{aligned}
$$
shows that $\phi_s$ is an algebra map.  A similar computation,
using that $\phi_{s-1}$ and $\theta_s$ commute with $d_X$, and
also that $d(V') = 0$, shows that $\phi_s$ also commutes with
$d_X$, and hence is a DG algebra map.  Thus, each $\phi_1, \dots,
\phi_r$ is a DG algebra map and in particular so is $\phi =
\phi_r$.

Next, we show the following:  That $\phi(v_1) = v_1 + v'_1$ and,
for $i = 2, \dots, r$,
$$\phi(v_i) = v_i + v'_i + I(v'_1, \dots, v'_{i-1}).
$$
This we do by induction on $s$.  Suppose inductively that we have
$\phi_s(v_1) = v_1 + v'_1$ and
$$\phi_s(v_i) = \begin{cases} v_i + v'_i + I(v'_1, \dots, v'_{i-1}) & \text{if $i = 2, \dots, s$} \\
v_i +  I(v'_1, \dots, v'_{i-1}) & \text{if $i = s+1, \dots, r$}
\end{cases}
$$
Induction starts with $s=1$, where the formulas
$$\phi_1(v_1) = v_1 + v'_1 \qquad \text{and} \qquad \phi_1(v_i) = v_i +
v'_1\theta_1(v_i)$$
give the result.  For the inductive step, we compute as follows:
$\phi_{s+1}(v_1) = \phi_s(v_1) + v'_{s+1}\theta_{s+1}(v_1) = v_1 +
v'_1$, since $1 < s+1$ and hence $\theta_{s+1}(v_1) = 0$.  For $i
= 2, \dots, s$, we have
$$\begin{aligned}
\phi_{s+1}(v_i) &= \phi_{s+1}(v_i) +
v'_{s+1}\theta_{s+1}\big(\phi_s(v_i)\big)\\
&= v_i + v'_i + I(v'_1, \dots, v'_{i-1}) +
v'_{s+1}\theta_{s+1}\big(v_i + v'_i + I(v'_1, \dots,
v'_{i-1})\big)\\
&= v_i + v'_i + I(v'_1, \dots, v'_{i-1})
\end{aligned}
$$
since $i < s+1$ and thus $\theta_{s+1}(v_i)= 0$, and also the
ideal $I(v'_1, \dots, v'_{i-1})$ is $\theta_{s+1}$-stable, as
$\theta_{s+1}(v'_i) = 0$.  Further, $\phi_{s+1}(v_{s+1}) =
\phi_{s+1}(v_{s+1}) +
v'_{s+1}\theta_{s+1}\big(\phi_s(v_{s+1})\big) = v_{s+1} +  I(v'_1,
\dots, v'_{s}) + v'_{s+1}\theta_{s+1}\big(v_{s+1} +  I(v'_1,
\dots, v'_{s})\big) = v_{s+1} + v'_{s+1} + I(v'_1, \dots,
v'_{s})$. Finally, for $i = s+2, \dots, r$, we have
$$\begin{aligned}
\phi_{s+1}(v_i) &= \phi_{s}(v_i) +
v'_{s+1}\theta_{s+1}\big(\phi_s(v_i)\big)\\
&= v_i + I(v'_1, \dots, v'_{i-1}) + v'_{s+1}\theta_{s+1}\big(v_i +
I(v'_1, \dots, v'_{i-1})\big)\\
&= v_i + I(v'_1, \dots, v'_{i-1})
\end{aligned}
$$
since $s+1 \leq i-1$.  This completes the induction.

Finally, we observe that, for any $z \in Z$, we have $\phi(z) \in
I(Z)$. This follows easily from the fact that $Z$ is
$\theta_i$-stable for each $i$.

From these facts, it is evident that $(1\cdot\epsilon)\circ\phi$
satisfies $(1\cdot\epsilon)\circ\phi(v_1) = v'_1$, and
$(1\cdot\epsilon)\circ\phi(v_i) = v'_i + I(v'_1, \dots, v'_{i-1})$
for $i = 2, \dots, r$.  It follows that
$(1\cdot\epsilon)\circ\phi$ is surjective.  Furthermore, we have
$(1\cdot\epsilon)\circ\phi(z) = 0$.  For the other projection, it
is evident from the definition of $\phi$ that we have $(\epsilon
\cdot1)\circ\phi = 1$.
\end{proof}

We deduce the following technical proposition.

\begin{proposition}\label{prop: no Lambda V cycles}
Suppose $V$ decomposes as $V = V'\oplus V''$, with $d_X(V') = 0$
and $V''$ satisfying the following:  For any cycle of the form $v
+ z + \chi$, with $v \in V$, $z \in Z$, and $\chi \in \land^{\geq
2}(V \oplus Z)$, we have $v \in V'$.  Suppose the complement $Z$
has been chosen to satisfy $\theta_i(Z) \subseteq \land V\otimes
\land^{+}Z$ for each $i$. Then any cycle of $\land^+(V \oplus Z)$
is in the ideal $I(V', Z)$ generated by $V'\oplus Z$.
\end{proposition}

\begin{proof}
The proof is similar to that of part (1) of \propref{prop:I(Z)
d-stable}. We argue by contradiction.  Suppose this is not true,
and that amongst cycles of the form $\alpha + \beta$, with $\alpha
\not= 0 \in \land V''$, $\beta \in I(V', Z)$, that the shortest
length term in any such $\alpha$ is $m \geq 2$.  Since each
$\theta_i$ commutes with the differential, $\theta_i(v')$ is a
cycle for each $i$.  Therefore, we must have that $\theta_i(V')
\subseteq \land^{\geq m} V'' + I(V', Z)$.  Now adjust our notation
slightly for this situation.  Write $V'' = \langle v''_1, \dots,
v''_s \rangle$ for suitable $s \leq r$, with corresponding
derivations $\theta''_1, \dots, \theta''_s$.  Let $\chi$ be a
cycle that displays a shortest length part in $\land V''$, and
suppose that $t \leq s$ is the highest index for which $v''_t$
occurs in this shortest length part.  Then write
$$\chi = \alpha' + \alpha'' v''_t + \alpha''' + \beta,$$
with $\alpha' \in \land^{m}(v''_1, \dots, v''_{t-1})$,
$\alpha''\not= 0 \in \land^{m-1}(v''_1, \dots, v''_{t-1})$,
$\alpha''' \in \land^{m+1}V''$, and $\beta \in I(V', Z)$. Since
$\theta''_t$ commutes with the differential, $\theta''_t(\chi)$ is
again a cycle.  However, we have $\theta''_t(\chi) = \alpha'' +
\theta''_t(\alpha''' + \beta)$ (recall that $\theta_i(v_j) = 0$
for $i > j$).  Using \lemref{lem:Z theta stable} and the fact that
$\theta_t$ is a derivation, we have $\theta_t(\alpha''' + \beta)
\in I(\land^{m}V'', V', Z)$.  This contradicts our minimal length
assumption.
\end{proof}

The next result is a consequence of Oprea's \thmref{thm:Oprea}.
In order to be self-contained we include here a short proof.

\begin{proposition}\label{prop: Gottlieb factor}
Suppose $\mathcal{M}_X$ is written as $\land(V'\oplus V'' \oplus
Z)$ as in \propref{prop: no Lambda V cycles}. Then we may identify
$V'$ with $\im h_X\circ (\omega_\Q)_\#$.  Furthermore, $X_\Q$
decomposes as a product $X_\Q \simeq S \times Y$ with $S$ a
product of odd-dimensional rational spheres whose minimal model is
$(\land V',0)$.
\end{proposition}

\begin{proof}
Suppose $(\land V,d)$ is a minimal model for $X$, $x\in V$ is a
cocycle of odd degree and that there is a derivation $\theta$ of
$\land V$ such that $[\theta, d] = 0$ and $\theta (x) = 1$. Write
$\land V = \land (x) \otimes \land W$. Then by induction on the
degree we can modify the choice of $W$ in order to have $d(W)
\subset \land W$, as follows. Suppose that $d(W^{<n}) \subset
\land W$ and let $v\in W^n$. We write $d(v) = x\alpha + \beta$
with $\alpha$ and $\beta \in \land W$. Then $[\theta , d]= 0$
implies $\alpha =-d(\theta (v))$. Replacing $v$ by $v' = v+
x\theta (v)$ we obtain $d(v') \in \land W$. In this way, we may
assume that $(\land V,d) = (\land x,0)\otimes (\land W,d)$, i.e.,
$X_\Q \simeq S^n_\Q \times Y$. Since we then have $G_*(X_\Q) \cong
G_*(S^n_\Q) \oplus G_*(Y)$, we can proceed in the same way with
$Y$. This results in the required decomposition.
\end{proof}

We may now prove \thmref{thm:mono} of the introduction.  In her
thesis \cite{ghorbal}, Sonia Ghorbal has obtained a criterion for
a map to be a homotopy monomorphism in the nilpotent category. In
order to be self-contained we reproduce here the statement and the
proof of this criterion.

\begin{proposition}[S.~Ghorbal]\label{prop:criterion}
Let $f \colon X \to Y$ be a map of rational spaces that admits a
minimal model of the form $\gamma \colon \big(\land(V\oplus
W),d\big)\to (\land V,\bar d)$ such that $\gamma(W) = 0$,
$\gamma(v) = v$ for $v\in V$, $d(W) \subseteq \land V \otimes
\land^{\geq 2}W$, and $d(V) \subseteq \land V + \land V \otimes
\land^{\geq 2}W$. Then $f$ is a homotopy monomorphism in the
nilpotent category.
\end{proposition}

\begin{proof}
We first recall from \cite{LMM} that two morphisms $k,l \colon
(\land V,d) \to (A,d)$ are homotopic if there exists a map $H \colon
\land(V\oplus\bar V\oplus V') \to (A,d)$ with $k(v) = H(v)$ and
$$l(v) = H(e^{sd+ds}(v))= H(v) +dH(\bar v) + \sum_{r\geq
1}\frac{1}{r!} H((sd)^r(v))$$
for each $v \in V$. In this definition $\bar Z$ and $Z'$ are graded
vector spaces, $\bar Z ^p = Z^{p+1}$, $(Z')^p = Z^p$. The
differential $d$ on $\land V$ is extended to $\land(V\oplus\bar
V\oplus V')$ by $d(\bar v) = v'$ and $d(v') = 0$ for $v \in V$. The
map $s$ is a degree $-1$ derivation defined by $s(v) = \bar v$ and
$s(\bar v) = s(v') = 0$.

Now suppose that $g,h \colon (\land V,\bar d)\to (\land T,d)$ are DG
algebra maps and that
$$\Phi \colon (\land (V\oplus W \oplus \bar V \oplus \bar W \oplus V'
\oplus W'),d) \to (\land T,d)$$
is a homotopy between $g\circ \gamma$ and $h\circ \gamma$. We denote
by $I$ the ideal of $\land(V\oplus W \oplus \bar V \oplus \bar W
\oplus V' \oplus W')$ generated by the vector spaces $\land^{2}W$,
$s(\land^{2}W)$, and $W'$.

First we show that $sd(I) \subseteq I$.  For this, write a typical
element of $I$ as $aA + bB + cC$ with $a \in \land^{2}W$, $b \in
s(\land^{2}W)$, $c \in W'$, and $A, B, C$ general elements of
$\land(V\oplus W \oplus \bar V \oplus \bar W \oplus V' \oplus W')$.
Then we have
$$sd(aA) = s\big((da)A \pm a (dA)\big) = (sda)A \pm (da)(sA) \pm (sa)(dA) \pm
a(sdA).$$
The last two terms are automatically in $I$.  The second, and hence
the first, is in $I$ due to the hypothesis on $d$.  A similar
analysis of the terms that occur shows that $sd(bB)$ and $sd(cC)$
are also in $I$.

Next, we show by induction that $\Phi(I) = 0$.  For this, choose a
basis $y_1, y_2, \ldots y_r,\ldots $ of $W$ with $\vert y_i\vert
\leq \vert y_{i+1}\vert$.  Also, denote by $W_{(n)}$ the subspace
$\langle y_1, \ldots, y_n \rangle$ of $W$ and by $I_{(n)}$ the
ideal generated by the vector spaces $\land^{2}(W_{(n)})$,
$s(\land^{2}(W_{(n)}))$, and $W_{(n)}'$.

When $n = 1$, we have $d(y_1) = 0$ from our hypothesis on $d$. Thus
we have
$$0 = h\circ\gamma(y_1) = \Phi(e^{sd+ds}(y_1)) = \Phi (y_1) + \Phi (y_1') =
g\circ\gamma(y_1) + \Phi(y_1') = \Phi(y_1')\,,$$
which starts the induction.  Now suppose that the result is true
for $i<n$. Then we have
$$0 = \Phi (e^{sd+ds}(y_n)) = \Phi (y_n) +
\Phi(y_n') + \sum_{r\geq 1} \frac{1}{r!}\Phi ((sd)^r(y_n))\,.$$
The hypothesis on $d$ implies that $sd(y_n) \in I_{(n-1)}$. A
refinement of the argument in the previous part shows that, in fact,
each $I_{(n-1)}$ is stable under $sd$.  Therefore, we have
$(sd)^r(y_n) \in I_{(n-1)}$ for $r\geq 1$. Since $\Phi(I_{(n-1)}) =
0$ by our induction hypothesis, we have that $\Phi((sd)^r(y_n)) = 0$
and therefore $\Phi(y_n') = 0$. Of course, $\Phi$ is already zero on
$W$ and hence vanishes on both $\land^{2}(W_{(n)})$ and
$s(\land^{2}(W_{(n)}))$.  Therefore, we have $\Phi(I_{(n)}) = 0$ and
the induction is complete.  It follows that the ideal $I$ possesses
two key properties, namely $sd(I) \subset I$ and $\Phi(I) = 0$. We
now define a homotopy
$$\Psi \colon (\land(V\oplus\bar V\oplus V'),\bar d) \to (\land
T,d)$$
simply by restricting $\Phi$. We remark that $(sd)^r(v)-(s\bar
d)^r(v) \in I$ for $v \in V$, for $r \geq 1$. Therefore the homotopy
ends at $\Psi (e^{s\bar d+\bar ds}(v))= \Phi (e^{sd+ds}(v)) =h(v)$.
Furthermore, we have $\Psi(v) = \Phi(v) = g(v)$ for $v\in V$. Thus
$\Psi$ is a homotopy between $g$ and $h$.

The argument so far shows that $f$ is a homotopy monomorphism in
the rational category. That is, if $A$ is any rational space, then
$$f_*\colon [A, X] \to [A, Y]$$
is one-to-one.  From the universal properties of localization, it
follows that $f$ is a homotopy monomorphism in the nilpotent
category.
\end{proof}

\begin{proof}[Proof of \thmref{thm:mono}]
For the ordinary evaluation map $\omega \colon \map(X,X;1) \to X$,
we have that $\Gamma_X \colon S_X \to X_\Q$ is a homotopy
monomorphism in the nilpotent category by \propref{prop:criterion}
and \propref{prop:I(Z) d-stable}. Now suppose that $\ev\colon E
\to X$ is any evaluation map. From \thmref{thm:Gamma+ factor}, we
have the following commutative diagram of solid arrows
$$\xymatrix{ & E\ar[r]^-{g} \ar[dd]_{\ev_\Q} \ar@{-->}@/_1pc/[dl]_{r_\ev}& \map(X,X_\Q;e)
\ar@{-->}@/^1pc/[dr]^-{r_X}\ar[dd]^{\omega_\Q}  \\
S_\ev \ar[dr]^{\Gamma_\ev}\ar[ur]_{\widetilde{\Gamma}_\ev}& & & S_X
\ar[dl]^{\Gamma_X}\ar[ul]^{\widetilde{\Gamma}_X}\\
 & X_\Q \ar@{=}[r] & X_\Q }$$
with retractions $r_X$ and $r_\ev$ of $\widetilde{\Gamma}_X$ and
$\widetilde{\Gamma}_\ev$ respectively. We define $j \colon S_\ev
\to S_X$ by $j = r_X\circ g\circ\widetilde{\Gamma}_\ev$ and claim
that this map admits a retraction. Recall that both $S_\ev$ and
$S_X$ are (finite) products of odd-dimensional rational spheres.
Also, since $\Gamma_\ev$ and $\Gamma_X$ are both injective in
rational homotopy and $\Gamma_X\circ j = \Gamma_\ev$, it follows
that $j$ is injective in rational homotopy.  In terms of minimal
models, then, we have a map $\mathcal{M}_j \colon (\land V, d=0)
\to (\land W, d=0)$ with $Q(\mathcal{M}_j)$ surjective.  But if
$Q(\mathcal{M}_j)$ is surjective, so too is $\mathcal{M}_j$.
Therefore, we may choose a splitting of $\mathcal{M}_j$ which
corresponds to a retraction of $j$.  Since $j$ admits a
retraction, it is a homotopy monomorphism. Finally, it follows
that $\Gamma_\ev$ is a composition of homotopy monomorphisms and
hence is a homotopy monomorphism.
\end{proof}

\begin{corollary}\label{cor:epi-mono}
Let $\ev\colon E \to X$ be any evaluation map.  Then $\ev_\Q$
factors as a composition $\ev_\Q = \Gamma_\ev\circ r_\ev$ with
$r_\ev$ a homotopy epimorphism and $\Gamma_\ev$ a homotopy
monomorphism in the nilpotent category.
\end{corollary}

\begin{proof}
The discussion at the start of this section concluded that $r_\ev$
is a homotopy epimorphism and the remainder follows immediately
from \thmref{thm:Gamma+ factor} and \thmref{thm:mono}.
\end{proof}

We remark that the fact that $\Gamma_\ev$ is associated to an
evaluation map is key in \thmref{thm:mono}.  In particular, we may
give the following example of a map $\gamma \colon S \to X$ from
an $H_0$-space $S$ into $X$ that is injective in rational homotopy
but is not a homotopy monomorphism in the nilpotent category.

\begin{example}
Let $S = S^3_a \times S^5$ and $X = S^3_a \vee S^3_b
\cup_{\alpha}e^8$, where $\alpha$ is the triple Whitehead bracket
$[a,[a,b]]$.  Then $\gamma\colon S \to X$ is an extension of
$(1\mid [a,b]) \colon S^3_a \vee S^5 \to X$ obtained using the
fact that $[a,[a,b]] = 0$ in $\pi_*(X)$.  Consider two maps $h, k
\colon  S^2 \times S^3 \to S^3_a\times S^5$. The map $h$ is the
composition
$$\xymatrix{S^2 \times S^3 \ar[r]^-{p_2} &  S^3 \ar[r]^-{i_1} &
S^3_a \times S^5}$$
and $k$ is the composition of the inclusion $S^3 \vee S^5 \to S^3
\times S^5$ with the map that consists of collapsing the cell
$S^2$ into a point:
$$\xymatrix{S^2 \times S^3 \ar[r]& S^2 \times S^3 / S^2 =
S^3 \vee S^5 \ar[r] & S^3 \times S^5.}$$
Clearly $h_\Q$ and $k_\Q$ are not homotopic because they do not
induce the same map in rational homology. However a simple
computation using minimal models show that the compositions
$f_\Q\circ h_\Q$ and $f_\Q\circ k_\Q$ are homotopic.
\end{example}

We finish this section with the topic of cyclic maps. A
\emph{cyclic map} $f \colon A \to X$ may be defined as a map that
lifts through the evaluation map $\omega \colon \map(X,X;1) \to
X$. This definition is easily seen to be equivalent to that given
above \thmref{thm:cyclic} \emph{via} the adjoint correspondence
between a map $A \to \map(X,X;1)$ that lifts $f$ and a map $A
\times X \to X$ that extends $(f\mid 1)$. Together with Sam Smith,
the second-named author has studied cyclic maps from the rational
homotopy point of view in \cite{Cyclic}. As we mentioned in the
introduction, our interest in the results of this paper arose from
that earlier work.

To state \corref{cor:G ev zero} we defined the Gottlieb groups of
a space relative to an evaluation map.     We say that a map $f
\colon A \to X$ is \emph{cyclic with respect to an evaluation map}
$\ev\colon E \to X$ if $f$ lifts through the evaluation map $\ev$.
Denote the set of homotopy classes of such maps by $G^\ev(A,X)$.
Upon rationalizing such a map, we obtain a map in
$G^{\ev_\Q}(A,X_\Q)$.

\begin{theorem}\label{thm:cyclic general}
Let $\ev\colon E \to X$ be an evaluation map with $X$ a nilpotent,
finite complex and let $A$ be a nilpotent space.  Then there are
bijections of sets
$$G^{\ev_\Q}(A,X_\Q) \cong [A, S_\ev] \cong \oplus_r \Hom
\big(H_r(A;\Q), G_r^\ev(X_\Q)\big)\,.$$
\end{theorem}

\begin{proof}
The first bijection is given by $(\Gamma_\ev)_* \colon [A, S_\ev]
\to G^{\ev_\Q}(A,X_\Q)$. This is a bijection by \thmref{thm:Gamma+
factor} and \thmref{thm:mono}. Now remark that $S_\ev$ has the
homotopy type of a product of rational \EM spaces, $S_\ev =
\prod_{i=1}^r K(\mathbb Q, n_i)$. By taking cohomology classes we
thus obtain a bijection
$$[A, S_\ev]
\stackrel{\cong}{\longrightarrow} \displaystyle\oplus_{i=1}^r
H^{n_i}(A;\mathbb Q)$$
and the result follows.
\end{proof}

Thus, for instance, we retrieve \cite[Th.3.2]{Cyclic}: If $A$ is a
space with non-zero rational cohomology in even degrees only, then
any map $g \colon A \to S_\ev$ must be null-homotopic, as $S_\ev$ is
a product of odd-dimensional rational \EM spaces. Consequently, this
hypothesis on $A$ entails the triviality of the set
$G^{\ev_\Q}(A,X_\Q)$. Many of the other results of \cite{Cyclic} may
be placed in context with the results of this paper.

If $X$ is a suspension, or more generally a co-$H_0$-space, then its
rationalized Gottlieb groups are generally trivial. Indeed, this is
the case as long as $X$ does not have the rational homotopy type of
a single sphere.  Therefore, it follows from \thmref{thm:w-sharp0}
that any cyclic map into a co-$H_0$-space that does not have the
rational homotopy type of a sphere is rationally trivial.  Basic
finiteness results, such as those of \cite{Lim}, follow from this.

Note, however, that a general cyclic map does not factor through the
product of odd spheres that corresponds to its image in rational
homotopy. That is, we are not able to extend \thmref{thm:Gamma+
factor} to cyclic maps. In particular, we note that there exist
cyclic maps that are trivial in rational homotopy and yet not
null-homotopic (e.g.~\cite[Ex.4.1]{Cyclic}).

\section{Evaluation Maps and Homology}\label{sec:cohomology}

After the preparatory results of \secref{sec:technical}, we prove in
this section the results concerning the homomorphism induced in
rational homology by an evaluation map.

\begin{proof}[Proof of \thmref{th:fact.homology}]
Consider $\omega \colon \map(X,X;1) \to  X$ as a special case
first. If $X$ is an $H_0$-space, then the multiplication of $X_\Q$
provides a section of $\omega_\Q$, so that $H_*(\omega;\Q)$ is
surjective.  If we have $X_\Q \simeq S^{2n+1}_\Q \times Y$, then
we may apply \thmref{thm: X = A x B}.  As $S^{2n+1}$ is an
$H_0$-space, the above observation gives that $\omega_{S^{2n+1}}$
is surjective on rational homology. Furthermore, the map $(i_1)^*$
in diagram (\ref{eq:w product}) immediately preceding \thmref{thm:
X = A x B} admits a section, namely $(p_1)^*$, and so it too is
surjective on rational homology.  It follows that $\im H_*(\omega
; \Q)$ contains at least the $H_*(S^{2n+1} ; \Q)$ factor and thus
is non-zero.  This establishes item (3) of
\thmref{th:fact.homology}.

Next, suppose that $h_X\circ (\omega_\Q)_\# = 0$.  We deduce from
  \lemref{lem:Z theta stable} and \propref{prop: no Lambda V cycles}
   that a model of $\Gamma_X$ is given by $$\mu : (\land
(V\oplus Z),d_X) \to (\land V,0)$$ with all cocycles of $\land
(V\oplus Z)$ in the ideal generated by $Z$ and $\mu(Z) = 0$.   Now
Proposition 3.2 (2) shows that the total Gottlieb element
$\Gamma_X$ induces the trivial homomorphism in rational
cohomology.

On the other hand, suppose that $h_X\circ (\omega_\Q)_\#$ has
image of dimension $r > 0$. Then \propref{prop: Gottlieb factor}
implies that we have $X_\Q \simeq S \times Y$ where $S$ is an
$r$-fold product of rational spheres of odd dimensions that
correspond to the image of $h_X\circ (\omega_\Q)_\#$.  Now we
apply \thmref{thm: X = A x B}   and conclude that $\im H_*(\omega
; \Q)$ contains   the $H_*(S; \Q)$ factor. Furthermore, we have
$h_Y\circ (\omega_Y)_\# = 0$, otherwise the image of $h_X\circ
(\omega_\Q)_\#$ would be of dimension $>r$.  Therefore,
$\widetilde{H}_*(\omega_Y; \Q) = 0$ and the image of
$H_*(\omega_\Q; \Q)$ is precisely the $H_*(S; \Q)$ factor. This
establishes the remaining items of \thmref{th:fact.homology} for
$\omega$.

Now consider a generalized evaluation map $\ev \colon E \to X$. We
suppose that $\im h_X\circ (\omega_\Q)_\#$ is of dimension $r$ and
$\im h_X\circ (\ev_\Q)_\#$ is of dimension $s$.  Since $\ev$
factors through $\omega$, we have $s\leq r$.  We write $X_\Q
\simeq S\times Y$ as above, and  we obtain a commutative diagram
$$\xymatrix{E \ar[r]^-{g} \ar[d]_{\ev_\Q} & \map(X,X_\Q; e)
\ar[d]^{\omega_\Q} \\
X_\Q \ar[r]_-{h}^-{\simeq} & S \times Y}$$
where $g$ is the $H$-map obtained from from the definition of a
generalized evaluation map. By Theorem 2.5    the coordinate maps
$p_1\circ\omega_\Q$ and $p_2\circ\omega_\Q$ factor through
$(\omega_S)_\Q$ and $(\omega_Y)_\Q$ respectively. Because of this
factorization, and  the fact that $\widetilde{H}_*(\omega_Y;\Q) =
0$, we may make the following identifications:
$$\im H_*(\ev_\Q;\Q) \cong \im H_*(\omega_\Q\circ g ;\Q)
\cong \im H_*(p_1\circ\omega_\Q\circ g;\Q) \subseteq H_*(S;\Q).$$
Since the composition $p_1\circ\omega_\Q\circ g\colon E \to S$
satisfies the hypotheses of \corref{cor:composition of H-map and
surj}, it admits a minimal model of the form $\varphi \colon
(\land V,0) \to (\land W,0)$ with $\varphi (V) \subset W$.   Then
the image of $p_1\circ\omega_\Q\circ g\colon E \to S$ in rational
homotopy has dimension $s$ and we may factor its minimal model
$\varphi \colon (\land V,0) \to (\land W,0)$ as the composition of
a surjection and an injection $\land(V_s\oplus K) \to \land V_s
\to \land(V_s\oplus K')$, with $V_s$ a vector space of dimension
$s$ isomorphic to the image of $\im h_X\circ (\ev_\Q)_\#$. This
corresponds to a factorization of $p_1\circ\omega_\Q\circ g\colon
E \to S$ as
$$\xymatrix{E\ar[rr]^-{p_1\circ\omega_\Q\circ g}\ar[dr]_q & &
S \simeq S'\times S''\\& S'\ar[ur]_{i_1}}$$
with $S'$ a product of odd-dimensional rational spheres with
minimal model $(\land V_s, 0)$.    It is now clear that the image
in homology of $\ev_\Q$ is isomorphic to $H_*(S';\Q)$.
\end{proof}

For $X$ a finite complex, a result of Gottlieb (\cite[Th.3]{Gott72})
says that if $\chi(X) \not=0$, then the first degree in which the
homomorphism induced by the evaluation map on rational cohomology
may be non-zero is even. With \thmref{th:fact.homology}, we sharpen
this result in a very significant way.

\begin{corollary}[\corref{cor:Gottlieb}]\label{cor:Gottlieb general}
Let $X$ be a nilpotent, finite space.  Suppose that $\chi(X) \not =
0$ or, more generally, that $X$ does not factor up to rational
homotopy as $X_\Q \simeq S^{2n+1}_\Q \times Y$. Then for every
evaluation map $\ev \colon E \to X$, we have
$\widetilde{H}_*(\ev;\Q)) = 0$.
\end{corollary}

Recall that $X$ is called a \emph{$c$-symplectic space} if it is an
even-dimensional rational Poincar{\'e} duality space that possesses
some class $x \in H^2(X;\Q)$, some power of which is a fundamental
class \cite{Lu-Op}.

\begin{corollary}[\corref{cor:symplectic}]\label{cor:symplectic general}
Let $X$ be a simply connected, $c$-symplectic space.  Then every
evaluation map $\ev \colon E \to X$ satisfies
$\widetilde{H}_*(\ev;\Q) = 0$.
\end{corollary}

\begin{proof}
It is evident that the cohomology algebra structure does not allow a
decomposition of the form $X \simeq_\Q S^{2n+1}\times X'$, and so
\thmref{th:fact.homology} implies the evaluation map is trivial in
rational homology.
\end{proof}

At the other extreme from the situation described in these
corollaries, we have the following:

\begin{corollary}\label{cor:H*(w) injective}
Let $\ev\colon E \to X$ be an evaluation map with $X$ a nilpotent,
finite complex. The following are equivalent:
\begin{enumerate}
\item The homomorphism $H_*(\ev) \colon H_*(E;\Q) \to
H_*(X;\Q)$ is surjective;
\item $\Gamma_\ev \colon S_\ev \to X$ is a rational homotopy
equivalence.
\end{enumerate}
When (1) and (2) pertain, $X$ is an $H_0$-space and the evaluation
map admits a section.
\end{corollary}

\begin{proof}
All parts follow easily from \thmref{thm:Gamma+ factor} and
\thmref{th:fact.homology}.
\end{proof}

\section{Conclusion: Some Open Problems}\label{sec:examples}

At present, we have very little information about the map $\ev
\colon \top(X,X;1) \to X$ or the other variations on the
evaluation map $\omega$ mentioned at the start of the
introduction.  It would be most interesting to identify
$G^\ev_*(X_\Q)$, the image in rational homotopy of $\ev$, or, more
generally the rational homotopy groups of $\top(X,X;1)$.  As
specific instances of this kind of problem, we offer the
following.

\begin{problem}
Let $M$ be a compact smooth manifold.  Is the image in (rational)
homotopy of the evaluation map $\ev\colon \diff(M,M;1) \to M$
strictly contained in, or equal to, the (rational) Gottlieb groups
of $M$?
\end{problem}

\begin{problem}
Let $X$ be an $H$-space and recall that $G_*(X) = \pi_*(X)$ in
this case. Let $H(X,X;1)$ denote the subspace of $\map(X,X;1)$
that consists of $H$-equivalences.  Is the evaluation map
$\ev\colon H(X,X;1) \to X$ surjective in (rational) homotopy?
\end{problem}

Assuming that $G^\ev_*(X)$ and $G_*(X)$ are generally different
from each other, it would be interesting to know whether there are
structural results for $G^\ev_*(X_\Q)$ comparable to those of
F{\'e}lix-Halperin for the ordinary Gottlieb groups.

\corref{cor:partial trivial} and \corref{cor:G ev zero} may be
used to give necessary conditions for certain maps to be the
connecting map of a fibration (cf.~\exref{ex:not partial}).  This
suggests the following particular version of an old problem of
Massey:

\begin{problem}
Let $p\colon \Omega B \to X$ be a map from a loop space to a
nilpotent, finite complex $X$.  When is $p$ the connecting map of
some fibration sequence $X \to E \to B$?
\end{problem}

It would be nice to find other situations in which the image in
rational homotopy groups of a map led to factorizations analogous
to those of \secref{sec:w-factor}.  In this direction, we offer
the following rather general problem:

\begin{problem}
Suppose given a map $f\colon X \to Y$ with $Y$ finite-dimensional.
If the image of $f_\#$ in rational homotopy groups is
finite-dimensional, does $f$ factor through an elliptic space?
\end{problem}

We have restricted ourselves entirely to the rational homotopy
context in this paper.  But it could be feasible to investigate
similar results working either integrally or localized at
different sets of primes.  We end with two ``moonshots" that
indicate how little we know outside the rational situation.

\begin{problem}
Let $X$ be a space with trivial Gottlieb groups (integrally).  Is
the evaluation map $\omega\colon \map(X,X;1) \to X$ null-homotopic?
\end{problem}

\begin{problem} Let $X$ be a nilpotent, finite complex.  When is a Gottlieb
element $S^n \to X$ a homotopy monomorphism, and not just a
rational homotopy monomorphism ?
\end{problem}

%\bibliographystyle{amsplain}
%\bibliography{evaluation}

\begin{thebibliography}{10}

\bibitem{F-H}
Y.~F{\'e}lix and S.~Halperin, \emph{Rational {LS} category and its
  applications}, Trans. Amer. Math. Soc. \textbf{273} (1982), no.~1, 1--38.
  \MR{84h:55011}

\bibitem{F-H-T}
Y.~F{\'e}lix, S.~Halperin, and J.-C. Thomas, \emph{Rational homotopy
theory},
  Graduate Texts in Mathematics, vol. 205, Springer-Verlag, New York, 2001.
  \MR{2002d:55014}

\bibitem{Gan1}
T.~Ganea, \emph{On monomorphisms in homotopy theory}, Topology
\textbf{6}
  (1967), 149--152. \MR{34 \#8402}

\bibitem{ghorbal}
S.~Ghorbal, \emph{Monomorphismes et {\'e}pimorphismes homotopiques},
  Ph.~D.~Thesis, Louvain-La-Neuve, 1996.

\bibitem{Got1}
D.~H. Gottlieb, \emph{A certain subgroup of the fundamental group},
Amer. J.
  Math. \textbf{87} (1965), 840--856.

\bibitem{Got2}
\bysame, \emph{On fibre spaces and the evaluation map}, Ann. Math.
\textbf{87}
  (1968), 42--55.

\bibitem{Go1}
\bysame, \emph{Evaluation subgroups of homotopy groups}, Amer. J.
Math.
  \textbf{91} (1969), 729--756.

\bibitem{Go2}
\bysame, \emph{Applications of bundle map theory},
Trans.~Amer.~Math.~Soc.
  \textbf{171} (1972), 23--50.

\bibitem{Gott72}
\bysame, \emph{The evaluation map and homology}, Michigan Math. J.
\textbf{19}
  (1972), 289--297. \MR{49 \#8005}

\bibitem{LMM}
S.~Halperin, \emph{Lectures on minimal models}, M\'em. Soc. Math.
France (N.S.)
  (1983), no.~9-10, 261. \MR{85i:55009}

 \bibitem{Hal2}
 \bysame, \emph{Torsion gaps in the homotopy of finite complexes},
 Topology
  \textbf{27} (1988), no.~3, 367--375. \MR{89h:55024}

\bibitem{H-M-R}
P.~Hilton, G.~Mislin, and J.~Roitberg, \emph{Localization of
nilpotent groups
  and spaces}, North-Holland Publishing Co., Amsterdam, 1975, North-Holland
  Mathematics Studies, No. 15, Notas de Matem\'atica, No. 55. [Notes on
  Mathematics, No. 55]. \MR{57 \#17635}

\bibitem{Lan}
G.~E. Lang, \emph{Localizations and evaluation subgroups}, Proc.
Amer. Math.
  Soc. \textbf{50} (1975), 489--494. \MR{51 \#4228}

\bibitem{Lim}
K.~L. Lim, \emph{On cyclic maps}, J. Austral. Math. Soc. Ser. A
\textbf{32}
  (1982), no.~3, 349--357. \MR{83e:55003}

\bibitem{Lu-Op}
G.~Lupton and J.~Oprea, \emph{Cohomologically symplectic spaces:
toral actions
  and the {G}ottlieb group}, Trans. Amer. Math. Soc. \textbf{347} (1995),
  no.~1, 261--288. \MR{95f:57056}

\bibitem{Cyclic}
G.~Lupton and S.~B. Smith, \emph{Cyclic maps in rational homotopy
theory},
  Math. Z. \textbf{249} (2005), no.~1, 113--124.

\bibitem{Mil}
J.~Milnor, \emph{On spaces having the homotopy type of {${\rm
CW}$}-complex},
  Trans. Amer. Math. Soc. \textbf{90} (1959), 272--280. \MR{20 \#6700}

\bibitem{Op86}
John Oprea, \emph{Decomposition theorems in rational homotopy
theory}, Proc.
  Amer. Math. Soc. \textbf{96} (1986), no.~3, 505--512. \MR{87h:55008}

\bibitem{Op87}
\bysame, \emph{The {S}amelson space of a fibration}, Michigan Math.
J.
  \textbf{34} (1987), no.~1, 127--141. \MR{88c:55015}

\bibitem{Var}
K.~Varadarajan, \emph{Generalised {G}ottlieb groups}, J.~Indian
Math. Soc.
  \textbf{33} (1969), 141--164.

\end{thebibliography}

%\end{document}
\providecommand{\bysame}{\leavevmode\hbox
to3em{\hrulefill}\thinspace}
\providecommand{\MR}{\relax\ifhmode\unskip\space\fi MR }
\providecommand{\MRhref}[2]{%
  \href{http://www.ams.org/mathscinet-getitem?mr=#1}{#2}
} \providecommand{\href}[2]{#2}

\end{document}